\documentclass{amsart}

\usepackage[leqno]{amsmath}
\usepackage{latexsym,amsthm,amsfonts,amssymb,xypic,rotating}
\usepackage{tikz}
\usetikzlibrary{matrix,arrows}
\usepackage{hyperref}

\newcommand{\numberseries}{\bfseries}   
\newlength{\thmtopspace}                
\newlength{\thmbotspace}                
\newlength{\thmheadspace}               
\newlength{\thmindent}                  
\newtheoremstyle{bfupright head,slanted body}
                {\thmtopspace}{\thmbotspace}
                {\slshape}{\thmindent}{\bfseries}{.}{\thmheadspace}
                {{\numberseries \thmnumber{#2\;}}\thmnote{#3}}
\newtheoremstyle{bfupright head,upright body}
                {\thmtopspace}{\thmbotspace}
                {\upshape}{\thmindent}{\bfseries}{.}{\thmheadspace}
                {{\numberseries \thmnumber{#2\;}}\thmnote{#3}}
\newtheoremstyle{fixed bf head,slanted body}
                {\thmtopspace}{\thmbotspace}{\slshape}
                {\thmindent}{\bfseries}{.}{\thmheadspace}
                {{\numberseries \thmnumber{#2\;}}\thmname{#1}\thmnote{ (#3)}}
\newtheoremstyle{fixed bf head,upright body}
                {\thmtopspace}{\thmbotspace}{\upshape}
                {\thmindent}{\bfseries}{.}{\thmheadspace}
                {{\numberseries \thmnumber{#2\;}}\thmname{#1}\thmnote{ (#3)}}
\newtheoremstyle{numbered paragraph}
                {\thmtopspace}{\thmbotspace}{\upshape}
                {\thmindent}{\upshape}{}{\thmheadspace}
                {{\numberseries \thmnumber{#2.}}}
\theoremstyle{bfupright head,slanted body}
\newtheorem{res}{}[section]             \newtheorem*{res*}{}
\theoremstyle{bfupright head,upright body}
\newtheorem{bfhpg}[res]{}               \newtheorem*{bfhpg*}{}
\theoremstyle{fixed bf head,slanted body}
\newtheorem{thm}[res]{Theorem}          \newtheorem*{thm*}{Theorem}
\newtheorem{prp}[res]{Proposition}      \newtheorem*{prp*}{Proposition}
\newtheorem{cor}[res]{Corollary}        \newtheorem*{cor*}{Corollary}
\newtheorem{lem}[res]{Lemma}            \newtheorem*{lem*}{Lemma}
\theoremstyle{fixed bf head,upright body}
       \newtheorem*{dfn*}{Definition}
\newtheorem{rmk}[res]{Remark}           \newtheorem*{rmk*}{Remark}
\newtheorem{stp}[res]{Setup}            \newtheorem*{stp*}{Setup}
\theoremstyle{numbered paragraph}
\newtheorem{ipg}[res]{}

\newlength{\thmlistleft}        
\newlength{\thmlistright}       
\newlength{\thmlistpartopsep}   
\newlength{\thmlisttopsep}      
\newlength{\thmlistparsep}      
\newlength{\thmlistitemsep}     

\newcounter{eqc} 
\newenvironment{eqc}{\begin{list}{\upshape (\textit{\roman{eqc}})}%
    {\usecounter{eqc}%
      \setlength{\leftmargin}{\thmlistleft}%
      \setlength{\labelwidth}{\thmlistleft}%
      \setlength{\rightmargin}{\thmlistright}%
      \setlength{\partopsep}{\thmlistpartopsep}%
      \setlength{\topsep}{\thmlisttopsep}%
      \setlength{\parsep}{\thmlistparsep}%
      \setlength{\itemsep}{\thmlistitemsep}}}%
  {\end{list}}%
\newcommand{\eqclbl}[1]{{\upshape(\textit{#1})}}
\newcounter{prt}
\newenvironment{prt}{\begin{list}{\upshape (\alph{prt})}%
    {\usecounter{prt}%
      \setlength{\leftmargin}{\thmlistleft}%
      \setlength{\labelwidth}{\thmlistleft}%
      \setlength{\rightmargin}{\thmlistright}%
      \setlength{\partopsep}{\thmlistpartopsep}%
      \setlength{\topsep}{\thmlisttopsep}%
      \setlength{\parsep}{\thmlistparsep}%
      \setlength{\itemsep}{\thmlistitemsep}}}%
  {\end{list}}%
\newcounter{rqm}
\newenvironment{rqm}{\begin{list}{\upshape (\arabic{rqm})}%
    {\usecounter{rqm}%
      \setlength{\leftmargin}{\thmlistleft}%
      \setlength{\labelwidth}{\thmlistleft}%
      \setlength{\rightmargin}{\thmlistright}%
      \setlength{\partopsep}{\thmlistpartopsep}%
      \setlength{\topsep}{\thmlisttopsep}%
      \setlength{\parsep}{\thmlistparsep}%
      \setlength{\itemsep}{\thmlistitemsep}}}%
  {\end{list}}%

\newenvironment{prf*}[1][Proof]{%
  \begin{proof}[\bf #1]
    \setcounter{equation}{0}
    }
  {\end{proof}
}

\newcommand{\pgref}[1]{\ref{#1}}

\newcommand{\stpref}[2][Setup~]{#1\ref{stp:#2}}
\newcommand{\thmref}[2][Theorem~]{#1\pgref{thm:#2}}
\newcommand{\corref}[2][Corollary~]{#1\pgref{cor:#2}}
\newcommand{\prpref}[2][Proposition~]{#1\pgref{prp:#2}}
\newcommand{\lemref}[2][Lemma~]{#1\pgref{lem:#2}}
\newcommand{\rmkref}[2][Remark~]{#1\pgref{rmk:#2}}
\newcommand{\secref}[2][Section~]{#1\ref{sec:#2}}
\newcommand{\tabref}[2][Table~]{#1\ref{tab:#2}}
\renewcommand{\eqref}[1]{(\pgref{eq:#1})}
\newcommand{\thmcite}[2][?]{\cite[Thm.~#1]{#2}}
\newcommand{\corcite}[2][?]{\cite[Cor.~#1]{#2}}
\newcommand{\prpcite}[2][?]{\cite[Prop.~#1]{#2}}
\newcommand{\lemcite}[2][?]{\cite[Lem.~#1]{#2}}
\newcommand{\dfncite}[2][?]{\cite[Def.~#1]{#2}}
\def\urltilda{\kern -.15em\lower .7ex\hbox{\~{}}\kern .04em} 
\newcommand{\set}[2][\mspace{1mu}]{\{#1 #2 #1\}}
\newcommand{\setof}[3][\mspace{1mu}]{\{#1#2 \mid #3#1\}}
\newcommand{\NN}{\mathbb{N}}
\newcommand{\ZZ}{\mathbb{Z}}
\newcommand{\qtext}[1]{\quad\text{#1}\quad}
\newcommand{\qqtext}[1]{\qquad\text{#1}\qquad}
\newcommand{\qand}{\qtext{and}}
\newcommand{\qqand}{\qqtext{and}}
\newcommand{\deq}{\:=\:}
\newcommand{\dge}{\:\ge\:}
\newcommand{\dle}{\:\le\:}
\newcommand{\gra}{\alpha}
\newcommand{\grb}{\beta}
\newcommand{\grg}{\gamma}
\newcommand{\grf}{\varphi}
\newcommand{\mfm}{\mathfrak{m}}
\newcommand{\is}{\cong}
\renewcommand{\le}{\leqslant}
\renewcommand{\ge}{\geqslant}
\newcommand{\lra}{\longrightarrow}
\newcommand{\Rhat}{\widehat{R}}
\newcommand{\mapdef}[4][\rightarrow]{\nobreak{#2\colon #3 #1 #4}}
\newcommand{\dif}[2][]{{\partial}^{#2}_{#1}}
\renewcommand{\H}[2][]{\operatorname{H}_{#1}(#2)}
\newcommand{\dptR}{\operatorname{depth}R}
\newcommand{\edim}[1]{\operatorname{edim}#1}
\newcommand{\dpt}[2][R]{\operatorname{depth}_{#1}#2}
\newcommand{\rnk}[2][k]{\operatorname{rank}_{#1}#2}
\newcommand{\pd}[2][R]{\operatorname{pd}_{#1}#2}
\newcommand{\Hom}[3][R]{\operatorname{Hom}_{#1}(#2,#3)}
\newcommand{\Ext}[4][R]{\operatorname{Ext}_{#1}^{#2}(#3,#4)}
\newcommand{\tp}[3][R]{\nobreak{#2\otimes_{#1}#3}}
\newcommand{\Tor}[4][R]{\operatorname{Tor}^{#1}_{#2}(#3,#4)}

\numberwithin{equation}{res}
\numberwithin{table}{section}

\newcommand{\Kzl}[1]{\operatorname{K}_{\sbt}(#1)}
\newcommand{\clC}[1]{\mathbf{C}(#1)} 
\newcommand{\clT}{\mathbf{T}} 
\newcommand{\clB}{\mathbf{B}} 
\newcommand{\clG}[1]{\mathbf{G}(#1)}

\newcommand{\clH}[1]{\mathbf{H}(#1)}

\newcommand{\ee}{\varepsilon} 
\newcommand{\e}{\mathsf{e}} 
\newcommand{\f}{\mathsf{f}} 
\newcommand{\g}{\mathsf{g}} 
\newcommand{\A}{\mathsf{A}} 
\newcommand{\B}{\mathsf{B}} 
\newcommand{\E}{\mathsf{E}} 
 
\newcommand{\F}{\mathsf{F}} 
\newcommand{\G}{\mathsf{G}}

\newcommand{\mfI}{\mathfrak{I}}

\newcommand{\mfM}{\mathfrak{M}}

\newcommand{\mfA}{\mathfrak{A}}
\newcommand{\mfB}{\mathfrak{B}}
\newcommand{\mfX}{\mathfrak{X}}
\newcommand{\sbt}{{\scriptscriptstyle\bullet}}
\newcommand{\prmm}[1]{m_{#1}}
\newcommand{\prmn}[1]{n_{#1}}
\newcommand{\prmp}[1]{p_{#1}}
\newcommand{\prmq}[1]{q_{#1}}
\newcommand{\prmr}[1]{r_{#1}}

\begin{document}

\title{Linkage classes of grade 3 perfect ideals}

\author[L.\,W. Christensen]{Lars Winther Christensen}

\address{Texas Tech University, Lubbock, TX 79409, U.S.A.}

\email{lars.w.christensen@ttu.edu}

\urladdr{http://www.math.ttu.edu/\urltilda lchriste}

\author[O. Veliche]{Oana Veliche}

\address{Northeastern University, Boston, MA~02115, U.S.A.}

\email{o.veliche@northeastern.edu}

\urladdr{https://web.northeastern.edu/oveliche}

\author[J. Weyman]{Jerzy Weyman}

\address{University of Connecticut, Storrs, CT~06269, U.S.A.}
\email{jerzy.weyman@uconn.edu}

\urladdr{http://www.math.uconn.edu/\urltilda weyman}

\thanks{This work is part of a body of research that started during
  the authors' visit to MSRI in Spring 2013 and continued during a
  months-long visit by L.W.C.\ to Northeastern University; the
  hospitality of both institutions is acknowledged with
  gratitude. L.W.C.\ was partly supported by NSA grant H98230-14-0140
  and Simons Foundation collaboration grant 428308, and J.W.\ was
  partly supported by NSF DMS grants 1400740 and 1802067.}

\date{4 June 2019}



\begin{abstract}
  While every grade $2$ perfect ideal in a regular local ring is
  linked to a complete intersection ideal, it is known not to be the
  case for ideals of grade $3$. We soften the blow by proving that
  every grade $3$ perfect ideal in a regular local ring is linked to a
  complete intersection or a Golod ideal. Our proof is indebted to a
  homological classification of Cohen--Macaulay local rings of
  codimension $3$. That debt is swiftly repaid, as we use linkage to
  reveal some of the finer structures of this classification.
\end{abstract}

\maketitle

\thispagestyle{empty}

\section{Introduction}

\noindent
Let $R$ be a local ring with maximal ideal $\mfm$ and residue field
$k = R/\mfm$. The difference between the \emph{embedding dimension}
and \emph{depth} of $R$,
\begin{equation*}
  \edim{R} = \rnk{\mfm/\mfm^2} \qqand
  \dptR = \inf\setof{i\in\ZZ}{\Ext{i}{k}{R} \ne 0}\:,
\end{equation*}
i.e.\ the number $c=\edim{R} - \dptR$, is called the \emph{embedding
  codepth} of $R$. Rings with $c=0$ are regular, rings with $c=1$ are
hypersurfaces, and for rings with $c=2$ there are two possibilities:
complete intersection or Golod. For $c=3$ the field of possibilities
widens: Such rings can be Gorenstein and not complete intersection, or
they may not even belong to any of the classes mentioned thus far.
There is, nevertheless, a classification of local rings of embedding
codepth $3$. It is based on multiplicative structures in homology, and
the details are discussed below.

In the 1980s Weyman~\cite{JWm89} and Avramov, Kustin, and Miller
\cite{AKM-88} established a classification scheme which---though it is
discrete in the sense that it does not involve moduli---is subtle
enough to facilitate a proof of the rationality of Poincar\'e series
of codepth $3$ local rings, an open question at the time. For this
purpose it was not relevant to know if local rings of every class in
the scheme actually exist, and it later turned out that they do
not. In 2012 Avramov \cite{LLA12} returned to the classification and
tightened it; that is, he limited the range of possible classes.  This
was necessary to use the classification to answer the codepth 3 case
of a question ascribed to Huneke about growth in the minimal injective
resolution of a local ring. In the same paper \cite[Question
3.8]{LLA12}, Avramov formally raised the question about realizability:
Which classes of codepth 3 local rings do actually occur? This paper
deals with aspects of the realizability question that can be studied
by linkage theory and, in turn, provide new insight on linkage of
grade $3$ perfect ideals.

To discuss our contributions towards an answer to the realizability
question, we need to describe the classification in further detail.
The $\mfm$-adic completion $\Rhat$ of $R$ is by Cohen's Structure
Theorem a quotient of a regular local ring. To be precise, there is a
regular local ring $Q$ with maximal ideal $\mfM$ and an ideal
$\mfI \subseteq \mfM^2$ with $Q/\mfI \is \Rhat$. The
Auslander-Buchsbaum Formula thus yields
$$c = \edim{Q} - \dpt[]{\Rhat} = \dpt[]{Q} - \dpt[Q]{\Rhat} =
\pd[Q]{\Rhat}\:.$$ For $c=3$, Buchsbaum and Eisenbud \cite{DABDEs77}
show that the minimal free resolution
\begin{equation*}
  F_\sbt \ = \ 0 \lra F_3 \lra F_2 \lra F_1 \lra F_0
\end{equation*}
of $\Rhat$ over $Q$ has a structure of a commutative differential
graded algebra. The multiplicative structure on $F_\sbt$ induces a
graded-commutative algebra structure on
$\Tor[Q]{\sbt}{\Rhat}{k} = \H[\sbt]{\tp[Q]{F_\sbt}{k}}$; as a
$k$-algebra, $\Tor[Q]{\sbt}{\Rhat}{k}$ is isomorphic to the Koszul
homology algebra $\H{K^R}$, in particular it is independent of the
presentation of $\Rhat$ and the choice of $F_\sbt$. The multiplicative
structure on $\Tor[Q]{\sbt}{\Rhat}{k}$ is the basis for the
classification. In \cite{AKM-88} the possible multiplication tables
are explicitly described, up to isomorphism; they are labeled $\clB$,
$\clC{3}$, $\clG{r}$, $\clH{p,q}$, and $\clT$, where the parameters
$p$, $q$, and $r$ are non-negative integers bounded by functions of
the invariants $m= \rnk[Q]{F_1}$, called the \emph{first derivation}
of $R$, and $n = \rnk[Q]{F_3}$, called the \emph{type} of $R$. Thus,
for fixed $m$ and $n$ there are only finitely many possible
structures, and \cite[Question 3.8]{LLA12} asks which ones actually
occur. It is known that a complete answer will have to take into
account the Cohen--Macaulay defect of the ring; in this paper we only
consider Cohen--Macaulay rings, and our main result pertains to the
two-parameter family:

\begin{thm}
  \label{thm:A}
  Let $R$ be a Cohen--Macaulay local ring of codimension 3 and class
  $\clH{p,q}$. Let $m$ and $n$ denote the first derivation and the
  type of $R$. The inequalities
  \begin{equation*}
    \prmp{} \dle m-1 \qqand \prmq{} \dle n
  \end{equation*}
  hold, and the following conditions are equivalent
  \begin{equation*}
    (i) \ \ \prmp{} = n+1 \qquad 
    (ii) \ \ \prmq{} = m-2 \qquad 
    (iii) \ \ \prmp{} = m-1 \qand\ \prmq{} = n\:.
  \end{equation*}
  Otherwise, i.e.\ when these conditions are not satisfied, there are
  inequalities
  \begin{equation*}
    p \le n-1 \qqand q\le m-4
  \end{equation*}
  with
  \begin{equation*}
    p=n-1 \ \ \text{ only if } \ \ q \equiv_2 m-4 \qqand 
    q=m-4 \ \ \text{ only if } \ \ p \equiv_2 n-1\:.
  \end{equation*}
\end{thm}
\noindent
The first set of inequalities in this theorem can be read off
immediately from the multiplication table for
$\Tor[Q]{\sbt}{\Rhat}{k}$, see \eqref{efg}, and the equivalence of conditions
\eqclbl{i}--\eqclbl{iii} is known from \corcite[3.3]{LLA12}. Thus,
what is new is the ``otherwise'' statement. The inequalities
$p \le n-1$ and $q\le m-4$ are proved in \secref{I}, and the last
assertions about equalities and congruences are proved in \secref{II}.

\begin{table}\small
  \begin{tikzpicture}[x=1.2cm,y=.6cm]
    \draw[step=1.0,gray,very thin](0,0) grid (7,6);
    \draw[black,thick](0,0) rectangle (7,6); \draw (0.5,0.5) node
    {$\clH{0,0}$}; \draw (1.5,0.5) node {$\clH{1,0}$}; \draw (2.5,0.5)
    node {$\clH{2,0}$}; \draw (3.5,0.5) node {$\clH{3,0}$}; \draw
    (0.5,1.5) node {$\clH{0,1}$}; \draw (1.5,1.5) node {$\clH{1,1}$};
    \draw (2.5,1.5) node {$\clH{2,1}$}; \draw (3.5,1.5) node
    {$\clH{3,1}$}; \draw (4.5,1.5) node {$\clH{4,1}$}; \draw (0.5,2.5)
    node {$\clH{0,2}$}; \draw (1.5,2.5) node {$\clH{1,2}$}; \draw
    (2.5,2.5) node {$\clH{2,2}$}; \draw (3.5,2.5) node {$\clH{3,2}$};
    \draw (0.5,3.5) node {$\clH{0,3}$}; \draw (2.5,3.5) node
    {$\clH{2,3}$}; \draw (4.5,3.5) node {$\clH{4,3}$}; \draw (6.5,5.5)
    node {$\clH{6,5}$}; \draw[white](0,0) -- (0,-.3);
  \end{tikzpicture}
  \caption{\small For $m=7$ and $n=5$ the innate bounds $p \le 6$ and $q\le 5$ are optimal. \thmref{A} exposes a finer structure to the classification that precludes $25$ of the $42$ classes $\clH{p,q}$ within these bounds. \vspace{-2.5\baselineskip}}  \label{tab:1}
\end{table}
To parse the theorem it may be helpful to visualize an instance of
it. For $m=7$ and $n=5$ the parameter $p$ can take values between $0$
and $m-1=6$, and $q$ can take values between $0$ and $n=5$; thus the
classes fit naturally in a $7\times 6$ grid; see \tabref{1}. It
follows from the theorem that at most $17$ of the $42$ fields in the
grid are populated, and in the experiments that informed the statement
of \thmref{A} we encountered all of them; we discuss this point
further in \secref{BGT}.

\begin{equation*}
  \ast \ \ \ast \ \ \ast
\end{equation*}
If $R$ and, therefore, $\Rhat$ is Cohen--Macaulay of codimension $3$,
then the defining ideal $\mfI$ with $\Rhat = Q/\mfI$ is perfect of
grade $3$. We obtain \thmref{A} as a special case of
results---\thmref[Theorems~]{pnqm} and \thmref[]{pnqm-1}---about grade
$3$ perfect ideals in regular local rings.

Our main tool of investigation is linkage. It is known that every
grade $2$ perfect ideal in $Q$ is linked to an ideal $\mfI$ such that
the local ring $Q/\mfI$ is complete intersection; such ideals are
called \emph{licci}. Not every grade $3$ perfect ideal is licci, see
for example \prpcite[3.5]{CVW-2}, but as a special case of
\thmref{link} we obtain:

\begin{thm}
  \label{thm:B}
  Every grade $3$ perfect ideal in a regular local ring $Q$ is linked
  to a grade $3$ perfect ideal $\mfI \subseteq Q$ such that $Q/\mfI$
  is complete intersection or Golod.
\end{thm}

The paper is organized as follows. We recall the details of the
classification from \cite{AKM-88,JWm89} in \secref{notation}. The
motor of the paper is \secref{link}; it has four statements that track
relations between the multiplicative structures on the Tor algebras of
linked ideals. As a first application of these statements, \thmref{B}
is proved in \secref{link}; the proofs of the four statements
themselves are deferred to the Appendix. The applications to the
classification scheme begin in \secref{dev2} and continue in
\secref[Sections~]{I} and \secref[]{II}, which contain the proof of
\thmref{A}. In \secref{BGT} we provide a summary of the status of the
realizability question.

\section{Multiplicative structures in homology}
\label{sec:notation}

\noindent
Throughout this paper, $Q$ denotes a commutative noetherian local ring
with maximal ideal $\mfM$ and residue field $k=Q/\mfM$. For an ideal
$\mfI \subseteq Q$ with $\pd[Q]{Q/\mfI} = 3$, let $F_\sbt \to Q/\mfI$
be a minimal free resolution over $Q$ and set
\begin{equation*}
  \prmm{\mfI} =  \rnk[Q]{F_1} \qqand \prmn{\mfI} = \rnk[Q]{F_3}\:;
\end{equation*}
if $Q$ is regular, then these numbers are the first derivation and the
type of $Q/\mfI$. Notice that one has $\rnk[Q]{F_0}=1$, which forces
$\rnk[Q]{F_2} = \prmm{\mfI} + \prmn{\mfI} - 1$ as $F_\sbt$ has Euler
characteristic $0$.

\begin{ipg}
  \label{ms}
  By a result of Buchsbaum and Eisenbud \cite{DABDEs77} the resolution
  $F_\sbt$ has a structure of a commutative differential graded
  algebra. This structure is not unique, but the induced
  graded-commutative algebra structure on
  $\A_\sbt = \H[\sbt]{\tp[Q]{F_\sbt}{k}} = \Tor[Q]{\sbt}{Q/\mfI}{k}$
  is unique.  By \cite{AKM-88} there exist bases
  \begin{equation*}
    \e_1,\ldots,\e_{\prmm{\mfI}} \ \text{ for $\A_1$}\,,\quad
    \f_1,\ldots, \f_{\prmm{\mfI}+\prmn{\mfI}-1} \ \text{ for $\A_2$}\,,\qand
    \g_1,\ldots,\g_{\prmn{\mfI}} \ \text{ for $\A_3$}
  \end{equation*}
  such that the multiplication on $\A_\sbt$ is one of following:
  \begin{equation}
    \label{eq:efg}
    \begin{aligned}
      \textbf{C}(3): \quad & \e_1\e_2 = \f_3 \ \ \e_2\e_3 = \f_1 \ \
      \e_3\e_1 = \f_2
      & \ \e_i\f_i = \g_1 \ \text{ for } \ 1\le i \le 3\\[.5ex]
      \textbf{T}: \quad & \e_1\e_2 = \f_3 \ \ \e_2\e_3 = \f_1 \ \
      \e_3\e_1  = \f_2\\[.5ex]
      \textbf{B}: \quad & \e_1\e_2 = \f_3 \
      & \ \e_i\f_i = \g_1 \ \text{ for } \ 1\le i \le 2\\[.5ex]
      \clG{r}: \quad & [r\ge 2]
      & \e_i\f_i = \g_1 \ \text{ for } \ 1\le i \le r\\[.5ex]
      \clH{p,q}: \quad & \e_{p+1}\e_i = \f_i \ \text{ for } \ 1\le
      i\le p & \e_{p+1}\f_{p+j} = \g_j \ \text{ for } \ 1\le j\le q
    \end{aligned}
  \end{equation}
  Here it is understood that all products that are not mentioned---and
  not given by those mentioned and the rules of graded
  commutativity---are zero. We say that $\mfI$, or $Q/\mfI$, is of
  class $\clC{3}$ if the multiplication on $\A_\sbt$ is given by
  $\clC{3}$ in \eqref{efg}; similarly for $\clB,\clG{r},\clH{p,q},$
  and $\clT$.
\end{ipg}

\begin{ipg}
  To deal with the multiplicative structures on $\A_\sbt$ it is
  helpful to consider a few additional invariants; set
  \begin{equation*}
    \prmp{\mfI} = \rnk{\A_1\!\cdot\!\A_1}\,, \quad
    \prmq{\mfI} = \rnk[k]{\A_1\!\cdot\!\A_2}\,, \qand
    \prmr{\mfI} = \rnk[k]{\delta_2^{\A}}
  \end{equation*}
  where $\mapdef{\delta_2^{\A}}{\A_2}{\Hom[k]{\A_1}{\A_3}}$ is defined
  by $\delta_2^{\A}(\f)(\e) = \f\e$ for $\e \in \A_1$ and
  $\f \in \A_2$. Depending on the class of $\mfI$, the values of these
  invariants are
  \begin{equation}
    \label{eq:pqr}
    \begin{array}{r|ccc}
      \text{Class of $\mfI$} & \prmp{\mfI} & \prmq{\mfI} & \prmr{\mfI}\\
      \hline
      \clB & 1 &1 &2 \\
      \clC{3} & 3 &1 &3 \\
      \clG{r}\ [r\ge 2]& 0 &1 &r \\
      \clH{p,q} & p & q & q\\
      \clT & 3 &0 &0 \\
    \end{array}
  \end{equation}  
\end{ipg}

\begin{ipg}
  \label{ci}
  Recall that an ideal $\mfX \subseteq Q$ is called \emph{complete
    intersection of grade $g$} if it is generated by a regular
  sequence of length $g$. For such an ideal the minimal free
  resolution of $Q/\mfX$ over $Q$ is the Koszul complex on the regular
  sequence. In particular, a complete intersection ideal is perfect. A
  perfect ideal of grade $g$ is called \emph{almost complete
    intersection} if it is minimally generated by $g+1$ elements.

  For a grade $3$ perfect ideal $\mfI \subseteq Q$ one has
  $\prmm{\mfI} \ge 3$, and the next conditions are equivalent; for the
  equivalence of \eqclbl{ii} and \eqclbl{iii} see the remark after
  \dfncite[2.2]{AKM-88}.

  \begin{eqc}
  \item $\prmm{\mfI} \le 3$.
  \item $\mfI$ is of class $\clC{3}$.
  \item $\mfI$ is complete intersection.
  \item $\prmm{\mfI} = 3$ and $\prmn{\mfI} = 1$.
  \end{eqc}
\end{ipg}

\begin{ipg}
  \label{gor}
  Recall that an ideal $\mfX \subseteq Q$ is called \emph{Gorenstein
    of grade $g$} if it is perfect of grade $g$ with $\prmn{\mfX}=1$,
  in which case one has $\Ext[Q]{g}{Q/\mfX}{Q} \is Q/\mfX$.

  If $\mfI \subseteq Q$ is Gorenstein of grade $3$, then $\mfI$ is of
  class $\clC{3}$ or of class $\clG{\prmm{\mfI}}$ with odd
  $\prmm{\mfI}\ge 5$; see the remark after \dfncite[2.2]{AKM-88}.
\end{ipg}

\enlargethispage*{\baselineskip} We refer to \cite[1.4]{LLA12} for the
following facts and precise references to their origins.

\begin{ipg}
  \label{reg}
  Assume that $Q$ is regular, and let $\mfI \subseteq Q$ be a grade
  $3$ perfect ideal.
  \begin{prt}
  \item The ring $Q/\mfI$ is complete intersection if and only if
    $\mfI$ is of class $\clC{3}$.
  \item The ring $Q/\mfI$ is Gorenstein and not complete intersection
    if and only if $\mfI$ is of class $\clG{\prmm{\mfI}}$ with
    $\prmn{\mfI}=1$.
  \item The ring $Q/\mfI$ is Golod if and only if $\mfI$ is of class
    $\clH{0,0}$.
  \end{prt}
\end{ipg}

\section{Multiplication in Tor algebras of linked ideals}
\label{sec:link}

\noindent
Let $\mfA \subseteq Q$ be a grade $3$ perfect ideal.  Recall that an
ideal $\mfB \subseteq Q$ is said to be \emph{directly linked} to
$\mfA$ if there exists a complete intersection ideal
$\mfX \subseteq \mfA$ of grade $3$ with $\mfB = (\mfX:\mfA)$. The
ideal $\mfB$ is then also a perfect ideal of grade $3$ with
$\mfX \subseteq \mfB$, and one has $\mfA = (\mfX:\mfB)$; see Golod
\cite{ESG80}. In particular, being directly linked is a reflexive
relation.  An ideal $\mfB$ is said to be \emph{linked} to $\mfA$ if
there exists a sequence of ideals
$\mfA = \mfB_0, \mfB_1,\ldots,\mfB_n = \mfB$ such that $\mfB_{i+1}$ is
directly linked to $\mfB_i$ for each $i=0,\ldots ,n-1$. Evidently,
being linked is an equivalence relation; the equivalence class of
$\mfA$ under this relation is called the \emph{linkage class} of
$\mfA$.

\begin{prp}
  \label{prp:linkBGT}
  Let $\mfA \subseteq Q$ be a grade $3$ perfect ideal.
  \begin{prt}
  \item If $\,\mfA$ is of class $\mathbf{B}$, then it is directly
    linked to a grade $3$ perfect ideal $\mfB$ with
    \begin{equation*}
      \prmm{\mfB} \deq \prmn{\mfA} + 2\,, \quad \prmn{\mfB} \deq
      \prmm{\mfA}-3\,, \qand \prmp{\mfB} \dge 2\:.
    \end{equation*}
    Moreover, $\mfB$ is of class $\mathbf{H}$.

  \item If $\,\mfA$ is of class $\mathbf{G}$, then it is directly
    linked to a grade $3$ perfect ideal $\mfB$ with
    \begin{equation*}
      \prmm{\mfB} \deq \prmn{\mfA} + 3\,, \quad \prmn{\mfB} \deq
      \prmm{\mfA}-3\,, \qand \prmp{\mfB} \dge \min\set{\prmr{\mfA},3}\:.
    \end{equation*}
    Moreover, $\mfB$ is of class $\mathbf{H}$.

  \item If $\,\mfA$ is of class $\mathbf{T}$ with $\prmm{\mfA} \ge 5$,
    then it is directly linked to a grade $3$ perfect ideal $\mfB$
    with
    \begin{equation*}
      \prmm{\mfB} \deq \prmn{\mfA} + 3\,, \quad \prmn{\mfB} \deq
      \prmm{\mfA}-3\,, \qand \prmq{\mfB} \dge 2\:.
    \end{equation*}
    In particular, $\mfB$ is of class $\mathbf{H}$.

  \item If $\,\mfA$ is of class $\mathbf{T}$, then it is directly
    linked to a grade $3$ perfect ideal $\mfB$ with
    \begin{equation*}
      \prmm{\mfB} \deq \prmn{\mfA} + 2\,, \quad \prmn{\mfB} \deq
      \prmm{\mfA}-3\,, \quad \prmq{\mfB} \deq 1\,, \qand \prmr{\mfB} \ge 2\:.
    \end{equation*}
    In particular, $\mfB$ is of class $\clB$ or $\mathbf{G}$.

  \item If $\,\mfA$ is of class $\mathbf{T}$, then it is directly
    linked to a grade $3$ perfect ideal $\mfB$ with
    \begin{equation*}
      \prmm{\mfB} \deq \prmn{\mfA} \qand \prmn{\mfB} \deq
      \prmm{\mfA}-3\:.
    \end{equation*}
  \end{prt}
\end{prp}

\begin{prf*}
  See \pgref{prf:linkBGT}.
\end{prf*}

The next three propositions deal with rings of class $\mathbf{H}$ the
way \prpref{linkBGT} deals with rings of class $\clB$, $\mathbf{G}$,
and $\clT$.

\begin{prp}
  \label{prp:linkH0}
  Let $\mfA \subseteq Q$ be a grade $3$ perfect ideal.  If\, $\mfA$ is
  of class $\mathbf{H}$ with $\prmm{\mfA}-3 \ge \prmp{\mfA}$, then it
  is directly linked to a grade $3$ perfect ideal $\mfB$ with
  \begin{equation*}
    \prmm{\mfB} = \prmn{\mfA} + 3\,, \quad \prmn{\mfB} =
    \prmm{\mfA}-3\,, \quad \prmp{\mfB} \ge \prmq{\mfA}\,, 
    \qand \prmq{\mfB} \ge \prmp{\mfA}\:.
  \end{equation*}
  Moreover, the following assertions hold:
  \begin{prt}
  \item If $\prmp{\mfA} \ge 1$, then $\prmp{\mfB} = \prmq{\mfA}$.
  \item If $\prmp{\mfA} \ge 2$, then $\mfB$ is of class
    $\clH{\prmq{\mfA},\;\cdot\;}$.
  \item If $\prmq{\mfA} \ge 2$, then $\prmq{\mfB} = \prmp{\mfA}$.
  \item If $\prmq{\mfA} \ge 3$, then $\mfB$ is of class
    $\clH{\;\cdot\;,\prmp{\mfA}}$.
  \item If $\prmq{\mfA} = 1$, then $\mfB$ is of class $\clB$ or
    $\mathbf{H}$.
  \item If $\prmq{\mfA} = 1$ and $\prmp{\mfA} = 0$, then $\mfB$ is of
    class $\mathbf{H}$.
  \item If $\mfB$ is of class $\mathbf{G}$, then
    $\prmr{\mfB} \le \prmm{\mfB}-2$.
  \end{prt}
\end{prp}

\begin{prf*}
  See \pgref{prf:linkH0}.
\end{prf*}

\begin{prp}
  \label{prp:linkH1}
  Let $\mfA \subseteq Q$ be a grade $3$ perfect ideal.  If\, $\mfA$ is
  of class $\mathbf{H}$ with $\prmm{\mfA}-2 \ge \prmp{\mfA} \ge 1$,
  then it is directly linked to a grade $3$ perfect ideal $\mfB$ with
  \begin{equation*}
    \prmm{\mfB} = \prmn{\mfA} + 2\,, \quad \prmn{\mfB} =
    \prmm{\mfA}-3\,, \quad \prmp{\mfB} \ge \prmq{\mfA}\,, 
    \qand \prmq{\mfB} \ge \prmp{\mfA}-1\:.
  \end{equation*}
  Moreover, the following assertions hold:
  \begin{prt}
  \item If $\prmp{\mfA} \ge 2$, then $\prmp{\mfB} = \prmq{\mfA}$.
  \item If $\prmp{\mfA} \ge 3$, then $\mfB$ is of class
    $\clH{\prmq{\mfA},\;\cdot\;}.$
  \item If $\prmq{\mfA} \ge 2$, then $\prmq{\mfB} = \prmp{\mfA}-1$.
  \item If $\prmq{\mfA} \ge 3$, then $\mfB$ is of class
    $\clH{\;\cdot\;,\prmp{\mfA}-1}$.
  \item If $\prmq{\mfA} = 1$, then $\mfB$ is of class $\clB$ or
    $\mathbf{H}$.
  \item If $\prmq{\mfA} = 1 = \prmp{\mfA}$, then $\mfB$ is of class
    $\mathbf{H}$.
  \item If $\prmn{\mfA} = 2$ and $\prmp{\mfB} = 3$, then $\mfA$ is of
    class $\clH{1,2}$ and $\mfB$ is of class $\mathbf{T}$.
  \end{prt}
\end{prp}

\begin{prf*}
  See \pgref{prf:linkH1}.
\end{prf*}

\begin{prp}
  \label{prp:linkH2}
  Let $\mfA \subseteq Q$ be a grade $3$ perfect ideal.  If\, $\mfA$ is
  of class $\mathbf{H}$ with $\prmm{\mfA}-1 \ge \prmp{\mfA} \ge 2$,
  then it is directly linked to a grade $3$ perfect ideal $\mfB$ with
  \begin{equation*}
    \prmm{\mfB} = \prmn{\mfA} + 1\,, \quad \prmn{\mfB} =
    \prmm{\mfA}-3\,, \quad \prmp{\mfB} \ge \prmq{\mfA}\,, 
    \qand \prmq{\mfB} \ge \prmp{\mfA}-2\:.
  \end{equation*}
  If $\prmn{\mfA} = 2$, then $\prmq{\mfA} = 2$ and $\mfB$ is complete
  intersection. Moreover, if $\prmm{\mfA} \ge 5$ or
  $\prmn{\mfA} \ge 3$, then the following assertions hold:
  \begin{prt}
  \item If $\prmp{\mfA} \ge 3$, then $\prmp{\mfB} = \prmq{\mfA}$.
  \item If $\prmp{\mfA} \ge 4$, then $\mfB$ is of class
    $\clH{\prmq{\mfA},\;\cdot\;}.$
  \item If $\prmq{\mfA} \ge 2$, then $\prmq{\mfB} = \prmp{\mfA}-2$.
  \item If $\prmq{\mfA} \ge 3$, then $\mfB$ is of class
    $\clH{\;\cdot\;,\prmp{\mfA}-2}$.
  \item If $\prmq{\mfA} = 1$, then $\mfB$ is of class $\clB$ or
    $\mathbf{H}$.
  \item If $\prmq{\mfA} = 1$ and $\prmp{\mfA}=2$, then $\mfB$ is of
    class $\mathbf{H}$.
  \end{prt}
\end{prp}

\begin{prf*}
  See \pgref{prf:linkH2}.
\end{prf*}

\begin{ipg}
  \label{Betti}
  For a grade $3$ perfect ideal $\mfI \subseteq Q$ the quantity
  $\prmm{\mfI} + \prmn{\mfI}$ is a measure of the size of the minimal
  free resolution $F_\sbt$ of $Q/\mfI$ over $Q$.  Indeed, one has
  \begin{equation*}
    \sum_{i=0}^3 \rnk[Q]{F_i} = 2(\prmm{\mfI} + \prmn{\mfI})\:;
  \end{equation*}
  we refer to this number as the \emph{total Betti number} of
  $\mfI$. As one has $\prmm{\mfI} \ge 3$ and $\prmn{\mfI} \ge 1$ the
  least possible total Betti number is $8$ and attained if and only if
  $\mfI$ is complete intersection; see \pgref{ci}.
\end{ipg}

The next corollary records the observation, already used in
\cite{AKM-88}, that any grade $3$ perfect ideal $\mfI$ with
$\prmp{\mfI} > 0$ is linked to an ideal with lower total Betti number.

\begin{cor}
  \label{cor:link}
  Let $\mfA \subseteq Q$ be a grade $3$ perfect ideal not of class
  $\clC{3}$.  There exists a grade $3$ perfect ideal $\mfB$ that is
  directly linked to $\mfA$ and has
  \begin{equation*}
    \prmm{\mfB} + \prmn{\mfB} \le 
    \prmm{\mfA} + \prmn{\mfA} - \min\set{2,\prmp{\mfA}}\:.
  \end{equation*}
\end{cor}

\begin{prf*}
  Immediate from \prpref[Propositions~]{linkBGT}--\prpref[]{linkH2}.
\end{prf*}

\thmref{B} from the introduction is a special case of the next result;
see \pgref{reg}.

\begin{thm}
  \label{thm:link}
  Every grade 3 perfect ideal in $Q$ is linked to a grade $3$ perfect
  ideal of class $\clC{3}$ or $\clH{0,0}$.
\end{thm}

\begin{prf*}
  Let $\mfA \subseteq Q$ be a grade $3$ perfect ideal, and assume that
  $\mfA$ is not of class $\clC{3}$ or $\clH{0,0}$. If $\prmp{\mfA}=0$,
  then $\mfA$ is of class $\mathbf{G}$ or $\mathbf{H}$ per
  \eqref{pqr}, so $\prmq{\mfA} \ge 1$. It now follows from
  \prpref{linkBGT}(b) and \prpref{linkH0} that $\mfA$ is linked to a
  grade $3$ perfect ideal $\mfB$ with $\prmp{\mfB} \ge 1$ and
  $\prmm{\mfB} + \prmn{\mfB} = \prmm{\mfA} + \prmn{\mfA}$.  Thus it
  follows from \corref{link} that every grade $3$ perfect ideal that
  is not of class $\clC{3}$ or $\clH{0,0}$ can be linked to a grade
  $3$ perfect ideal with smaller total Betti number, and ideals of
  class $\clC{3}$ have the smallest possible total Betti number; see
  \pgref{Betti}.
\end{prf*}

\section{Grade 3 perfect ideals generated by at most 5 elements}
\label{sec:dev2}

\noindent
Under the assumption that $Q$ is Gorenstein, the next result can be
deduced from Avramov's proof of \thmcite[2]{LLA81a}; see also
\cite[3.4.2]{LLA12}. The Gorenstein assumption is used to invoke a
result of Buchsbaum and Eisenbud \cite{DABDEs77}, but it follows from
later work of Golod \cite{ESG80} it is superfluous; see also
Brown~\cite[Intro.\ to Sec.~2]{AEB87}. Here we give a proof based on
the results in \secref{link}.

\begin{thm}
  \label{thm:m4}
  Let $\mfA \subseteq Q$ be a grade 3 perfect ideal with
  $\prmm{\mfA} = 4$.
  \begin{prt}
  \item If $\prmn{\mfA}$ is odd, then $\prmn{\mfA} \ge 3$ and $\mfA$
    is of class $\clT$.
  \item If $\prmn{\mfA} \ge 4$ is even, then $\mfA$ is of class
    $\clH{3,0}$.
  \item If $\prmn{\mfA} =2$, then $\mfA$ is of class $\clH{3,2}$.
  \end{prt}
  In particular, one has $\prmp{\mfA}=3$ and
  $\prmq{\mfA}\in\set{0,2}$.
\end{thm}

\begin{prf*}
  By \pgref{ci} the ideal $\mfA$ is not of class $\clC{3}$. If $\mfA$
  is of class $\clB$ or $\mathbf{G}$, then there exists by
  \prpref{linkBGT}(a,b) a grade $3$ perfect ideal $\mfB$ of class
  $\mathbf{H}$ with $\prmn{\mfB} = 1$; that is, a Gorenstein ideal of
  class $\mathbf{H}$. Per \pgref{gor} no such ideal exists, so $\mfA$
  is of class $\mathbf{H}$ or $\mathbf{T}$, cf.~\eqref{efg}.

  (a): Assume that $\prmn{\mfA}$ is odd. It is immediate from
  \pgref{ci} and \pgref{gor} that $\prmn{\mfA}$ cannot be $1$, so
  $\prmn{\mfA} \ge 3$ holds. To prove that the ideal $\mfA$ is of
  class $\mathbf{T}$, assume towards a contradiction that it is of
  class $\mathbf{H}$. By \eqref{efg} one has
  $\prmp{\mfA} \le \prmm{\mfA} - 1 = 3$.  If $\prmp{\mfA} \le 1$
  holds, then there exists by \prpref{linkH0} a grade $3$ perfect
  ideal $\mfB$ with $\prmn{\mfB} = 1$ and
  $\prmm{\mfB}= \prmn{\mfA}+3$, which contradicts \pgref{gor} as
  $\prmn{\mfA}+3$ is even. If $\prmp{\mfA} \ge 2$ holds, then
  \prpref{linkH2} yields a similar contradiction.

  (b): Assume that $\prmn{\mfA}$ is even. To prove that the ideal
  $\mfA$ is of class $\mathbf{H}$, assume towards a contradiction that
  it is of class $\mathbf{T}$. It follows from \prpref{linkBGT}(d)
  that there exists a grade $3$ perfect ideal $\mfB$ with
  $\prmn{\mfB} = 1$ and $\prmm{\mfB}= \prmn{\mfA}+2$, which
  contradicts \pgref{gor} as $\prmn{\mfA}+2$ is even. Thus $\mfA$ is
  of class $\mathbf{H}$ with $\prmp{\mfA} \le 3$, see \eqref{efg}, and
  we argue that equality holds.  If $\prmp{\mfA} = 0$, then $\mfA$ is
  by \prpref{linkH0} linked to a grade $3$ perfect ideal $\mfB$ with
  $\prmn{\mfB} = 1$. By \pgref{gor} the ideal $\mfB$ is of class
  $\clG{\prmm{\mfB}}$ which contradicts \prpref[]{linkH0}(g).  If
  $1 \le \prmp{\mfA} \le 2$ holds, then there exists by
  \prpref{linkH1} a grade $3$ perfect ideal $\mfB$ with
  $\prmn{\mfB}=1$ and $\prmm{\mfB} = \prmn{\mfA} + 2$, which
  contradicts \pgref{gor} as $\prmn{\mfA}+2$ is even.  Thus, $\mfA$ is
  of class $\mathbf{H}$ with $\prmp{\mfA} = 3$.

  Assume now that $\prmn{\mfA} \ge 4$ holds. By \prpref{linkH2}(a) the
  ideal $\mfA$ is linked to a grade $3$ perfect ideal $\mfB$ with
  $\prmn{\mfB} = 1$ and $\prmp{\mfB} = \prmq{\mfA}$. By \pgref{gor}
  and \eqref{efg} one has $\prmp{\mfB}=0$, so $\mfA$ is of class
  $\clH{3,0}$.

  (c): The argument above shows that $\mfA$ is of class $\mathbf{H}$
  with $\prmp{\mfA} = 3$. As $\prmn{\mfA} = 2$ \prpref{linkH2} yields
  $\prmq{\mfA} = 2$.
\end{prf*}

The results in \secref{link} easily yield the codimension $3$ case of
the fact that every almost complete intersection ideal is linked to a
Gorenstein ideal; see \prpcite[5.2]{DABDEs77}.

\begin{rmk}
  \label{rmk:m4}
  Let $\mfA \subseteq Q$ be a grade $3$ perfect ideal with
  $\prmm{\mfA}=4$. It follows from \thmref{m4}, \prpref{linkBGT}(e),
  \prpref{linkH2}, and \pgref{gor} that $\mfA$ is directly linked to a
  Gorenstein ideal $\mfB$ with
  \begin{equation*}
    \prmm{\mfB} = 
    \begin{cases}
      \prmn{\mfA} & \text{ if $\prmn{\mfA}$ is odd;}\\
      \prmn{\mfA}+1 & \text{ if $\prmn{\mfA}$ is even.}
    \end{cases}
  \end{equation*}
\end{rmk}

From the proof of the theorem in J.~Watanabe's \cite{JWt73} one can
deduce that every almost complete intersection ideal is linked to a
complete intersection ideal. Here is an explicit statement.

\begin{prp}
  Let $\mfA \subseteq Q$ be a grade $3$ perfect ideal with
  $\prmm{\mfA}=4$. There exists a grade $3$ perfect ideal of class
  $\clC{3}$ that is linked to $\mfA$ in at most $\prmn{\mfA}-2$ links
  if $\prmn{\mfA}$ is odd and at most $\prmn{\mfA}-1$ links if
  $\prmn{\mfA}$ is even.
\end{prp}

\begin{prf*}
  First assume that $\prmn{\mfA}$ is even. If $\prmn{\mfA}=2$, then it
  follows from \rmkref{m4} that $\mfA$ is directly linked to a
  Gorenstein ideal $\mfB$ with $\prmm{\mfB} = 3$, i.e.\ $\mfB$ is of
  class $\clC{3}$; see \pgref{ci}. Now let $n\ge 2$ be an integer and
  assume that the statement holds for ideals $\mfA'$ with
  $\prmn{\mfA'}=2(n-1)$. If $\prmn{\mfA}=2n$, then it follows from
  \rmkref{m4} that $\mfA$ is directly linked to a Gorenstein ideal
  $\mfB$ with $\prmm{\mfB} = 2n + 1$.  By \prpref{linkBGT}(b) there is
  an ideal $\mfA'$ that is directly linked to $\mfB$ and has
  $\prmm{\mfA'} = 4$ and $\prmn{\mfA'} = 2(n-1)$. By assumption
  $\mfA'$ is linked to an ideal of class $\clC{3}$ in at most
  $2(n-1) - 1$ links, so $\mfA$ is linked to the same ideal in at most
  $2 + 2(n-1) - 1 = \prmn{\mfA} - 1$ links.

  Now assume that $\prmn{\mfA}$ is odd. If $\prmn{\mfA}=3$, then it
  follows from \rmkref{m4} that $\mfA$ is directly linked to a
  Gorenstein ideal $\mfB$ with $\prmm{\mfB} = 3$, i.e.\ $\mfB$ is of
  class $\clC{3}$; see \pgref{ci}. If $\prmn{\mfA}=2n+1$ for some
  $n\ge 2$, then it follows from \rmkref{m4} that $\mfA$ is directly
  linked to a Gorenstein ideal $\mfB$ with $\prmm{\mfB} = 2n + 1$.  By
  \prpref{linkBGT}(b) there is an ideal $\mfA'$ that is directly
  linked to $\mfB$ and has $\prmm{\mfA'} = 4$ and
  $\prmn{\mfA'} = 2(n-1)$. By what has already been proved, $\mfA'$ is
  linked to an ideal of class $\clC{3}$ in at most $2(n-1) - 1$ links,
  so $\mfA$ is linked to the same ideal in at most
  $2 + 2(n-1) - 1 = \prmn{\mfA} - 2$~links.
\end{prf*}

The next result is proved by Brown \thmcite[4.5]{AEB87}. For
completeness we include a proof based on the results in
\secref{link}.

\begin{prp}
  \label{prp:n2}
  Let $\mfA \subseteq Q$ be a grade 3 perfect ideal with
  $\prmm{\mfA} \ge 5$ and $\prmn{\mfA}=2$. If $\prmp{\mfA} \ge 1$,
  then the following assertions hold:
  \begin{prt}
  \item $\mfA$ is of class $\clB$ if $\prmm{\mfA}$ is odd.
  \item $\mfA$ is of class $\clH{1,2}$ if $\prmm{\mfA}$ is even.
  \end{prt}
  In particular, one has $\prmp{\mfA} = 1$.
\end{prp}

\begin{prf*}
  By the assumptions $\prmm{\mfA} \ge 5$ and $\prmn{\mfA}=2$, the ideal
  $\mfA$ is not of class $\clC{3}$, see \pgref{ci}, so the assumption
  $\prmp{\mfA} \ge 1$ implies that $\mfA$ is of class $\mathbf{B}$,
  $\mathbf{H}$, or $\mathbf{T}$; see \eqref{pqr}.

  If $\mfA$ is of class $\mathbf{T}$, then there exists by
  \prpref{linkBGT}(d) a grade $3$ perfect ideal $\mfB$ of class
  $\mathbf{B}$ or $\mathbf{G}$ with $\prmm{\mfB} = 4$, but by
  \thmref{m4} no such ideal exists.

  If $\mfA$ is of class $\mathbf{B}$, then there exists by
  \prpref{linkBGT}(a) a grade $3$ perfect ideal $\mfB$ of class
  $\mathbf{H}$ with $\prmm{\mfB} = 4$, and
  $\prmn{\mfB} = \prmm{\mfA} - 3$, so $\prmm{\mfA}$ is odd by
  \thmref{m4}.

  Assume now that $\mfA$ is of class $\mathbf{H}$. By \eqref{efg} one
  has $\prmp{\mfA} \le \prmm{\mfA} - 1$, so if $\prmp{\mfA} \ge 2$
  holds, then there exists by \prpref{linkH2} a grade $3$ perfect
  ideal $\mfB$ with $\prmm{\mfB} = 3$ and
  $\prmn{\mfB} = \prmm{\mfA} - 3 \ge 2$, which contradicts
  \pgref{ci}. Thus one has $\prmp{\mfA} = 1$, and it now follows from
  \prpref{linkH1} that $\mfA$ is linked to a grade $3$ perfect ideal
  $\mfB$ with $\prmm{\mfB} = 4$, and $\prmn{\mfB} = \prmm{\mfA} -
  3$. By \thmref{m4} one has $\prmp{\mfA} = 3$, so by
  \prpref[]{linkH1}(g) the ideal $\mfA$ is of class $\clH{1,2}$ and
  $\mfB$ is of class $\mathbf{T}$. Now it follows from \thmref{m4}
  that $\prmm{\mfA}$ is even. The in particular statement is now
  immediate, see \eqref{pqr}.
\end{prf*}

\prpref{linkBGT} applies to restrict the possible classes of five
generated ideals.

\begin{thm}
  \label{thm:m5}
  Let $\mfA \subseteq Q$ be a grade 3 perfect ideal with
  $\prmm{\mfA}=5$. The following assertions hold:
  \begin{prt}
  \item $\mfA$ is of class $\mathbf{B}$ only if $\prmn{\mfA} = 2$.
  \item $\mfA$ is of class $\mathbf{G}$ only if $\prmn{\mfA} = 1$.
  \item $\mfA$ is of class $\mathbf{T}$ only if $\prmn{\mfA} \ge 4$.
  \end{prt}
\end{thm}

Notice that by part (b) a five generated ideal of class $\mathbf{G}$
is Gorenstein and hence of class $\clG{5}$; see \pgref{gor}.

\begin{prf*}
  (a): If $\mfA$ is of class $\mathbf{B}$, then it is by
  \prpref{linkBGT}(a) linked to a grade $3$ perfect ideal $\mfB$ with
  \begin{equation*}
    \prmm{\mfB} \deq \prmn{\mfA} + 2\,, \quad \prmn{\mfB} \deq
    2\,, \qand \prmp{\mfB} \dge 2\:.
  \end{equation*}
  If $\prmn{\mfA} = 1$, then one has $\prmm{\mfB} = 3$ and, therefore,
  $\prmn{\mfB} = 1$ per \pgref{ci}; a contradiction.  If
  $\prmn{\mfA} \ge 3$, then one has
  $\prmm{\mfB} = \prmn{\mfA}+2 \ge 5$ so \prpref{n2} yields
  $\prmp{\mfB} = 1$; a contradiction.
  
  (b): If $\mfA$ is of class $\mathbf{G}$, then it is by
  \prpref{linkBGT}(b) linked to a grade $3$ perfect ideal $\mfB$ with
  \begin{equation*}
    \prmm{\mfB} \deq \prmn{\mfA} + 3\,, \quad \prmn{\mfB} \deq
    2\,, \qand \prmp{\mfB} \dge \min\set{\prmr{\mfA},3} \ge 2\:.
  \end{equation*}
  If $\prmn{\mfA} \ge 2$ holds, then one has $\prmm{\mfB} \ge 5$ and
  hence $\prmp{\mfB} = 1$ by \prpref{n2}; a contradiction.

  (c): If $\mfA$ is of class $\mathbf{T}$, then it is by
  \prpref{linkBGT}(e) linked to a grade $3$ perfect ideal $\mfB$ with
  $\prmm{\mfB} = \prmn{\mfA}$ and $\prmn{\mfB} = 2$; by \pgref{ci}
  this forces $\prmn{\mfA} \ge 4$.
\end{prf*}

The next result was first proved by S\'anchez \cite{RSn89}.

\begin{prp}
  \label{prp:sanchez}
  Let $\mfA \subseteq Q$ be a grade $3$ perfect ideal. If
  $\prmm{\mfA} \ge 5$ and $\prmn{\mfA}=3$, then $\prmp{\mfA} \ne
  3$. In particular, $\mfA$ is not of class $\clT$.
\end{prp}

\begin{prf*}
  Assume towards a contradiction that $\prmm{\mfA} \ge 5$ and
  $\prmn{\mfA} = 3 = \prmp{\mfA}$ hold. Per \eqref{pqr} the ideal
  $\mfA$ is of class $\mathbf{H}$ or $\clT$.

  If $\mfA$ is of class $\mathbf{H}$, then there exists by
  \prpref{linkH2}(a) a grade $3$ perfect ideal $\mfB$ with
  $\prmm{\mfB} \deq 4$ and $\prmp{\mfB} \deq \prmq{\mfA}$.  By
  \thmref{m4} one has then $\prmp{\mfB}=3$ and, therefore,
  $\prmq{\mfA}=3$. Thus \prpref[]{linkH2}(c) yields
  $\prmq{\mfB} = \prmp{\mfA} - 2 = 1$, which contradicts
  \thmref[]{m4}.

  If $\mfA$ is of class $\clT$, then it follows from
  \prpref{linkBGT}(d) there exists a grade $3$ perfect ideal $\mfB$ in
  $Q$ of class $\clB$ or $\mathbf{G}$ with
  \begin{equation*}
    \prmm{\mfB} = 5\,, \quad \prmn{\mfB} = \prmm{\mfA}-3\,, \quad \prmq{\mfB} = 1\,, 
    \qand \prmr{\mfB} \ge 2\:.
  \end{equation*}
  By \thmref{m5} one has $\prmm{\mfA} \ge 6$ and, therefore,
  $\prmn{\mfB} \ge 3$, which contradicts \thmref[]{m5}.
\end{prf*}

\section{Class H: Bounds on the parameters of multiplication}
\label{sec:I}

\noindent
Throughout this section, the local ring $Q$ is assumed to be regular.

\begin{lem}
  \label{lem:square}
  Let $\mfI$ be a grade $3$ perfect ideal in $Q$. If $\mfI$ is not
  contained in $\mfM^2$, then one has $\prmm{\mfI} - \prmn{\mfI} = 2$.
\end{lem}

\begin{prf*}
  Choose an element $x \in \mfI \setminus \mfM^2$ and set $P = Q/(x)$;
  notice that $x$ is a $Q$-regular element. By the
  Auslander--Buchsbaum Formula one has \mbox{$\pd[P]{Q/\mfI} = 2$}, so
  the format of the minimal free resolution of $Q/\mfI$ as a
  $P$-module is
  \begin{equation*}
    0 \lra P^{a-1} \lra P^a \lra P
  \end{equation*}
  for some $a\ge 2$. Now it follows from \thmcite[2.2.3]{ifr} that the
  minimal free resolution of $Q/\mfI$ as a $Q$-module has format
  $0 \to Q^{a-1} \to Q^{2a-1} \to Q^{a+1} \to Q$. In particular, one
  has $\prmm{\mfI} = a+1$ and $\prmn{\mfI} = a-1$.
\end{prf*}

\begin{ipg}
  Let $\mfI \subseteq Q$ be a grade $3$ perfect ideal. If $\mfI$ is
  not contained in $\mfM^2$, then $\prmn{\mfI}=1$ implies that $\mfI$
  is of class $\clC{3}$ by \lemref{square} and \pgref{ci}. If $\mfI$
  is contained in $\mfM^2$, then $\prmn{\mfI}=1$ holds if and only if
  $\mfI$ is Gorenstein, i.e.\ of class $\clC{3}$ or
  $\clG{\prmm{\mfI}}$; see~\pgref{gor}. In particular, one has the
  following lower bounds on $\prmm{\mfI}$ and $\prmn{\mfI}$ depending
  on the class of $\mfI$.
  \begin{equation}
    \label{eq:mn}
    \begin{array}{r|rr}
      \text{Class of $\mfI$}   & \prmm{\mfI} & \prmn{\mfI}\\
      \hline
      \clB    & \ge 4 & \ge 2\\
      \clC{3} & 3 & 1 \\
      \clG{r} & \ge 4 & \ge 1\\
      \clH{p,q}  & \ge 4 & \ge 2\\
      \clT  & \ge 4 & \ge 2\\
    \end{array}
  \end{equation}  
\end{ipg}

\lemref{square} together with results from \cite{LLA12} yield bounds
on the parameters of multiplication for ideals of class \textbf{H}.

\begin{prp}
  \label{prp:lucho}
  Let $\mfI \subseteq Q$ be a grade $3$ perfect ideal of class
  $\mathbf{H}$; one has
  \begin{alignat}{3}
    \label{eq:pmqn}
    \prmp{\mfI} &\dle \prmm{\mfI}-1 &\qqand \prmq{\mfI} &\dle \prmn{\mfI}\:;\\
    \label{eq:pnqm}
    \prmp{\mfI} &\dle \prmn{\mfI}+1 &\qqand \prmq{\mfI} &\dle
    \prmm{\mfI}-2\:.
  \end{alignat}
  Moreover, the following
  conditions are equivalent:
  \begin{equation*}
    (i) \ \ \prmp{\mfI} = \prmn{\mfI} + 1 \qquad 
    (ii) \ \ \prmq{\mfI} = \prmm{\mfI} - 2 \qquad 
    (iii) \ \ \prmp{\mfI} = \prmm{\mfI}-1 \qand\ \prmq{\mfI} = \prmn{\mfI}\:.
  \end{equation*}
  Notice that when these conditions are satisfied, one has
  $\prmm{\mfI} - \prmn{\mfI} = 2$.
\end{prp}

\begin{prf*}
  The inequalities \eqref{pmqn} are immediate from \eqref{efg}. If
  $\mfI$ is not contained in $\mfM^2$, then by \lemref{square} the
  inequalities in \eqref{pnqm} agree with those in \eqref{pmqn}. If
  $\mfI$ is contained in $\mfM^2$, then the inequalities \eqref{pnqm}
  hold by \thmcite[3.1]{LLA12}. Moreover, for such an ideal the
  conditions $(i)$--$(iii)$ are equivalent by \corcite[3.3]{LLA12}.

  Now let $\mfI$ be an ideal that is not contained in $\mfM^2$; in
  view of \lemref{square} it suffices to prove that conditions $(i)$
  and $(ii)$ are equivalent. By \eqref{mn} one has
  $\prmn{\mfI} \ge 2$, and we proceed by induction on
  $\prmn{\mfI}$. For $\prmn{\mfI} = 2$ one has $\prmm{\mfI} = 4$ by
  \lemref{square}, and the equivalence of $(i)$ and $(ii)$ is trivial,
  as $\mfI$ is of class $\clH{3,2}$ by \thmref{m4}. Now let
  $\prmn{\mfI} = n \ge 3$, by \lemref{square} one has
  $\prmm{\mfI} = n+2$.
  
  $(i)\!\!\implies\!\!(ii)$: Assume that $\prmp{\mfI} = n+1$ holds. By
  \prpref{linkH2}(b) there exists a grade $3$ perfect ideal $\mfB$ of
  class $\mathbf{H}$ with
  \begin{equation*}
    \prmm{\mfB} = n+1\,,\quad \prmn{\mfB} = n-1\,, \qand 
    \,\prmq{\mfB} \ge n-1\,.
  \end{equation*}
  By \eqref{pmqn} one has $\prmq{\mfB} \le n-1$, so
  $\prmq{\mfB} = n - 1 = \prmm{\mfB} - 2$ holds. By the induction
  hypothesis, or by the equivalence of $(i)$--$(iii)$ for ideals
  contained in $\mfM^2$, one now has
  $\prmp{\mfB} = \prmn{\mfB} + 1 = n$. As $\prmp{\mfB} = \prmq{\mfI}$
  holds by \prpref[]{linkH2}(a), one has $\prmq{\mfI} = n = \prmm{\mfI}-2$.

  $(ii)\!\!\implies\!\!(i)$: Assume that $\prmq{\mfI} = n$ holds. If
  one has $\prmp{\mfI} \le 1$, then there exists by \prpref{linkH0}(d)
  a grade $3$ perfect ideal $\mfB$ of class $\mathbf{H}$ with
  $\prmm{\mfB} = n+3$, $\prmn{\mfB} = n-1$, and $\prmp{\mfB} \ge
  n$. By \eqref{pnqm} one has $\prmp{\mfB} \le \prmn{\mfB} + 1 = n$,
  so equality holds. It follows from \lemref{square} that $\mfB$ is
  contained in $\mfM^2$, so conditions $(i)$--$(iii)$ are
  equivalent. However, as $\prmm{\mfB} -1 =n+2$ the equality
  $\prmp{\mfB} = \prmm{\mfB} -1$ does not hold; a contradiction. Thus,
  one has $\prmp{\mfI} \ge 2$, and it follows from \prpref{linkH2}(d)
  that there exists a grade $3$ perfect ideal $\mfB$ of class
  $\mathbf{H}$ with
  \begin{equation*}
    \prmm{\mfB} = n+1\,,\quad \prmn{\mfB} = n-1,\, \qand
    \,\prmp{\mfB} \ge n\,.
  \end{equation*}
  By \eqref{pmqn} one has $\prmp{\mfB} \le n$, so
  $\prmp{\mfB} = n = \prmn{\mfB} + 1$ holds. By the induction
  hypothesis, or by the equivalence of $(i)$--$(iii)$ for ideals
  contained in $\mfM^2$, one now has
  $\prmq{\mfB} = \prmm{\mfB} -2 = n-1$. As
  $\prmq{\mfB} = \prmp{\mfI} -2$ holds by \prpref[]{linkH2}(c), one has
  $\prmp{\mfI} = n+1$.
\end{prf*}

The goal of this section is to establish part of \thmref{A} by proving
the following statement:

\begin{thm}
  \label{thm:pnqm}
  Let $\mfI \subseteq Q$ be a grade $3$ perfect ideal of class
  $\mathbf{H}$. If $\prmp{\mfI} \ne \prmn{\mfI} + 1$ or
  $\prmq{\mfI} \ne \prmm{\mfI}-2$ hold, then one has
  $\prmp{\mfI} \le \prmn{\mfI} - 1$ and
  $\prmq{\mfI} \le \prmm{\mfI} - 4$.
\end{thm}

The proof takes up the balance of the section; it is propelled by
\prpref[Propositions~]{linkH0}, \prpref[]{linkH1}, and
\prpref[]{linkH2}. For example, an ideal $\mfA$ with $\prmm{\mfA} = 6$
is by these results linked to an ideal $\mfB$ with $\prmn{\mfB}=3$,
and $\prmp{\mfB}$ and $\prmq{\mfB}$ are essentially determined by
$\prmq{\mfA}$ and $\prmp{\mfA}$, respectively. Thus restrictions on
$\prmq{\mfB}$ (or $\prmp{\mfB)}$ imply restrictions on
$\prmp{\mfA}$ (or $\prmq{\mfA}$). Here is an outline: Let
$\mfI \subseteq Q$ be a grade $3$ perfect ideal of class
$\mathbf{H}$. By \prpref{lucho} one can assume that
$\prmp{\mfI} \ne \prmn{\mfI}+1$ and $\prmq{\mfI} \ne \prmm{\mfI}-2$
hold, and it suffices to prove
\begin{center}
  (1) \ \ $\prmp{\mfI} \ne \prmn{\mfI}$ \hspace{1.5pc} and
  \hspace{1.5pc} (2) \ \ $\prmq{\mfI} \ne \prmm{\mfI}-3\:.$
\end{center}
By \eqref{mn} one has $\prmm{\mfI} \ge 4$ and $\prmn{\mfI} \ge 2$. If
$\prmm{\mfI} = 4$ or $\prmn{\mfI} = 2$ holds, then both (1) and (2)
follow from \thmref{m4} and \prpref{n2}. A few other low values of
$\prmm{\mfI}$ and $\prmn{\mfI}$ require special attention, and after
that the proof proceeds by induction: (2) for $\prmn{\mfI} = 3$
implies (1) for $\prmm{\mfI} = 6$, which in turn implies (2) for
$\prmn{\mfI} = 4$ etc.

\begin{lem}
  \label{lem:m=5}
  For every grade 3 perfect ideal $\mfA \subseteq Q$ with
  $\prmm{\mfA}=5$, one has $\prmq{\mfA} \ne 2$.
\end{lem}

\begin{prf*}
  If $\mfA$ is not of class $\mathbf{H}$, then per \eqref{pqr} one has
  $\prmq{\mfA} \le 1$, so assume that $\mfA$ is of class $\mathbf{H}$.
  Towards a contradiction, assume that $\prmq{\mfA}=2$ holds. By
  \eqref{pmqn} one has
  \begin{equation*}
    \tag{1}
    \prmp{\mfA} \le 4 \qqand \prmn{\mfA} \ge 2\:.
  \end{equation*}

  If $\prmp{\mfA} \le 2$ holds, then as $\prmq{\mfA}=2$ by assumption
  it follows from \prpref{linkH0} that there exists a grade $3$
  perfect ideal $\mfB$ in $Q$ with
  \begin{equation*}
    \prmm{\mfB} = \prmn{\mfA}+3\,,\quad \prmn{\mfB} = 2\,, \qand \prmp{\mfB} \ge 2\:.
  \end{equation*}
  By $(1)$ one has $\prmm{\mfB} = \prmn{\mfA} + 3 \ge 5$, so
  \prpref{n2} yields $\prmp{\mfB} = 1$; a contradiction.

  If $\prmp{\mfA} = 3$ holds, then as $\prmq{\mfA}=2$ by assumption,
  \prpref{linkH1}(b,c) yields the existence a grade $3$ perfect ideal
  $\mfB$ in $Q$ of class $\clH{2,2}$ with $\prmn{\mfB} =2$.  As
  $\prmn{\mfB} \ne 1$ one has $\prmm{\mfB} \ge 4$, see \eqref{mn}, and
  it follows from \thmref{m4} and \prpref{n2} that no such ideal
  $\mfB$ exists.

  If $\prmp{\mfA}=4$ holds, then as $\prmq{\mfA}=2$ by assumption,
  \prpref{linkH2}(b,c) yields the existence of a grade $3$ perfect
  ideal $\mfB$ in $Q$ of class $\clH{2,2}$ with $\prmn{\mfB} =2$.  As
  $\prmn{\mfB} \ne 1$ one has $\prmm{\mfB} \ge 4$, see \eqref{mn}, and
  it follows from \thmref{m4} and \prpref[]{n2} that no such ideal
  $\mfB$ exists.
\end{prf*}

\begin{lem}
  \label{lem:n}
  Let $n\ge 3$ be an integer and assume that for every grade $3$
  perfect ideal $\mfB \subseteq Q$ with $\prmn{\mfB} = n$ one has
  $\prmp{\mfB} \ne \prmn{\mfB}$. For every grade 3 perfect ideal
  $\mfA \subseteq Q$ with $\prmm{\mfA} = n + 3$ one then has
  $\prmq{\mfA} \ne \prmm{\mfA} - 3$.
\end{lem}

\begin{prf*}
  Let $n\ge 3$ be an integer and $\mfA \subseteq Q$ a grade 3 perfect
  ideal with $\prmm{\mfA} = n+3$.  Assume towards a contradiction that
  $\prmq{\mfA} = \prmm{\mfA} - 3 = n$ holds. As $n\ge 3$, it follows
  that $\mfA$ is of class $\mathbf{H}$, see \eqref{pqr}. Per
  \eqref{pmqn} one has
  \begin{equation*}
    \tag{1}
    \prmp{\mfA} \le n+2 \qqand \prmn{\mfA} \ge n\:.
  \end{equation*}
  We treat the cases $\prmp{\mfA} \le n$, $\prmp{\mfA}=n+1$, and
  $\prmp{\mfA}=n+2$ separately.

  \emph{Case 1.} If $\prmp{\mfA} \le n = \prmm{\mfA}-3$ holds, then as
  $\prmq{\mfA} = n$ by assumption it follows from \prpref{linkH0}(d)
  that there exists a grade $3$ perfect ideal $\mfB$ in $Q$ of class
  $\mathbf{H}$ with
  \begin{equation*}
    \tag{2}
    \prmm{\mfB} = \prmn{\mfA}+3\,, \quad \prmn{\mfB} = n\,, \qand 
    \prmp{\mfB} \ge n\:.
  \end{equation*}
  From $(2)$ and $(1)$ one gets
  \begin{equation*}
    \prmm{\mfB} - \prmn{\mfB} = \prmn{\mfA}+3-n \ge 3\:,
  \end{equation*}
  so \prpref{lucho} yields
  $\prmp{\mfB} \le \prmn{\mfB}$.  By $(2)$ one now has
  $\prmp{\mfB} = \prmn{\mfB}$, which is a contradiction.

  \emph{Case 2.} If $\prmp{\mfA} = n+1 = \prmm{\mfA}-2$ holds, then as
  $\prmq{\mfA} = n$ by assumption it follows from \prpref{linkH1}(b,c)
  that there exists a grade $3$ perfect ideal $\mfB$ in $Q$ of class
  $\clH{n,n}$ with $\prmn{\mfB} = n$, and that is a contradiction.

  \emph{Case 3.} If $\prmp{\mfA} = n+2 = \prmm{\mfA}-1$ holds, then as
  $\prmq{\mfA} = n$ by assumption it follows from \prpref{linkH2}(b,c)
  that there exists a grade $3$ perfect ideal $\mfB$ in $Q$ of class
  $\clH{n,n}$ with $\prmn{\mfB} = n$, and that is a contradiction.
\end{prf*}

\begin{lem}
  \label{lem:m}
  Let $m\ge 6$ be an integer and assume that for every grade 3 perfect
  ideal $\mfA \subseteq Q$ with $m-1 \le \prmm{\mfA} \le m$ one has
  $\prmq{\mfA} \ne \prmm{\mfA} - 3$. For every grade $3$ perfect ideal
  $\mfB \subseteq Q$ with $\prmn{\mfB} = m - 2$ one then has
  $\prmp{\mfB} \ne \prmn{\mfB}$.
\end{lem}

\begin{prf*}
  Let $m\ge 6$ be an integer and $\mfB \subseteq Q$ a grade 3 perfect
  ideal with $\prmn{\mfB} = m - 2$.  Assume towards a contradiction
  that
  \begin{equation*}
    \tag{1}
    \prmp{\mfB} = \prmn{\mfB} = m-2
  \end{equation*}
  holds. As $m-2 \ge 4$, it follows that $\mfB$ is of class
  $\mathbf{H}$, see \eqref{pqr}, and \eqref{pmqn} yields
  \begin{equation*}
    \tag{2}
    m-2 = \prmp{\mfB} \le \prmm{\mfB} -1\:.
  \end{equation*}
  We treat the cases $\prmm{\mfB} = m$ and $\prmm{\mfB} \ne m$
  separately.

  \emph{Case 1.} Assume that $\prmm{\mfB} = m$ holds.  Per $(0)$ one
  has $\prmp{\mfB} = \prmm{\mfB} - 2$, so by \prpref{linkH1}(b) there
  exists a grade $3$ perfect ideal $\mfA$ in $Q$ of class $\mathbf{H}$
  with
  \begin{equation*}
    \prmm{\mfA} = m\,,\quad \prmn{\mfA} = m-3\,,\qand
    \prmq{\mfA} \ge m-3\:.
  \end{equation*}
  In view of \eqref{pnqm} one now has
  \begin{equation*}
    \prmm{\mfA} - 2 \ge \prmq{\mfA} \ge \prmm{\mfA} - 3\:.
  \end{equation*}
  As $\prmm{\mfA} = m$ one has $\prmq{\mfA} \ne \prmm{\mfA}-3$ by
  assumption, which forces $\prmq{\mfA} = \prmm{\mfA} - 2$.  As
  one has $\prmm{\mfA} - \prmn{\mfA} = 3$, this contradicts \prpref{lucho}.

  \emph{Case 2.} Assume that $\prmm{\mfB} \ne m$ holds. Per (1) one
  has $\prmp{\mfB} \le \prmm{\mfB}-1$, so by $(0)$ and
  \prpref{linkH2}(b) there exists a grade $3$ perfect ideal $\mfA$ in
  $Q$ of class $\mathbf{H}$ with
  \begin{equation*}
    \prmm{\mfA} = m-1\,,\quad \prmn{\mfA} = \prmm{\mfB}-3\,\qand
    \prmq{\mfA} \ge m-4\:.
  \end{equation*}
  In view of \eqref{pnqm} one now has
  \begin{equation*}
    \prmm{\mfA} - 2 \ge \prmq{\mfA} \ge \prmm{\mfA} - 3\:.
  \end{equation*}
  As $\prmm{\mfA} = m-1$, one has $\prmq{\mfA} \ne \prmm{\mfA}-3$ by
  assumption, so that forces $\prmq{\mfA} = \prmm{\mfA} - 2$.
  As one has $\prmm{\mfA} - \prmn{\mfA} = m - \prmm{\mfB} + 2 \ne 2$,
  this contradicts.
\end{prf*}

\begin{prf*}[Proof of \ref{thm:pnqm}]
  By \eqref{pnqm} there are inequalities
  $\prmp{\mfI} \le \prmn{\mfI} + 1$ and
  $\prmq{\mfI} \le \prmm{\mfI}-2$. By the assumptions and
  \prpref{lucho} neither equality holds, so it is sufficient to show
  that $\prmp{\mfI} \ne \prmn{\mfI}$ and
  $\prmq{\mfI} \ne \prmm{\mfI}-3$ hold.  Per \eqref{mn} one has
  $\prmm{\mfI} \ge 4$ and $\prmn{\mfI} \ge 2$; to prove the assertion
  we show by induction on $l \in\NN$ that the following hold:
  \begin{prt}
  \item Every grade $3$ perfect ideal $\mfA \subseteq Q$ of class
    $\mathbf{H}$ with $\prmm{\mfA} = l+3$ has $\prmq{\mfA} \ne l$.
  \item Every grade $3$ perfect ideal $\mfB \subseteq Q$ of class
    $\mathbf{H}$ with $\prmn{\mfB} = l+1$ has $\prmp{\mfB} \ne l+1$.
  \end{prt}

  For $l=1$ part (a) follows from \thmref{m4} and (b) follows from
  \thmref{m4} and \prpref{n2}.

  For $l=2$ part (a) is \lemref{m=5}, and (b) holds by \thmref{m4} and
  \prpref{sanchez}.

  Let $l \ge 3$ and assume that (a) and (b) hold for all lower values
  of $l$. By part (b) for $l-1$ the hypothesis in \lemref{n} is
  satisfied for $n=l$, and it follows that (a) holds for $l$. By part
  (a) for $l$ and $l-1$ the hypothesis in \lemref{m} is satisfied for
  $m=l+3$, and it follows that (b) holds for $l$.
\end{prf*}

\begin{cor}
  \label{cor:MT1m2}
  Let $\mfI \subseteq Q$ be a grade $3$ perfect ideal of class
  $\mathbf{H}$. If $\prmm{\mfI} - \prmn{\mfI} \ne 2$ holds, then one
  has $\prmp{\mfI} \le \prmn{\mfI} - 1$ and
  $\prmq{\mfI} \le \prmm{\mfI} - 4$.
\end{cor}

\begin{prf*}
  The assertion is immediate from \prpref{lucho} and \thmref{pnqm}.
\end{prf*}

\section{Class H: Extremal values of the parameters of multiplication}
\label{sec:II}

\noindent
Throughout this section, the local ring $Q$ is assumed to be
regular. We establish the remaining part of \thmref{A} by proving the
following:

\begin{thm}
  \label{thm:pnqm-1}
  For every grade $3$ perfect ideal $\mfI \subseteq Q$ of class
  $\mathbf{H}$ the following hold
  \begin{prt}
  \item If $\prmp{\mfI} \ge \prmn{\mfI}-1$, then
    $\prmq{\mfI} \equiv_2 \prmm{\mfI}-4$.
  \item If $\prmq{\mfI} \ge \prmm{\mfI}-4$, then
    $\prmp{\mfI} \equiv_2 \prmn{\mfI}-1$.
  \end{prt}
\end{thm}

The proof follows the template from \secref{I}. It is, eventually, an
induction argument with the low values of $\prmm{\mfI}$ and
$\prmn{\mfI}$ handled separately in \thmref{m4}, \prpref{n2}, and
\lemref[Lemmas~]{m5}--\lemref[]{n4}.

\begin{lem}
  \label{lem:m5}
  For every grade $3$ perfect ideal $\mfA \subseteq Q$ of class
  $\mathbf{H}$ with $\prmm{\mfA} = 5$ and $\prmq{\mfA} \ge 1$ one has
  $\prmp{\mfA} \equiv_2 \prmn{\mfA}-1$.
\end{lem}

\begin{prf*}
  By \eqref{mn} and \eqref{pmqn} one has $\prmn{\mfA} \ge 2$ and
  $\prmp{\mfA} \le 4$. The assumption $\prmq{\mfA} \ge 1$ together
  with \eqref{pnqm} and \lemref{m=5} yield $\prmq{\mfA} = 1$ or
  $\prmq{\mfA} = 3$. We treat odd and even values of $\prmn{\mfA}$
  separately.

  Assume that $\prmn{\mfA}$ is even. The goal is to show that
  $\prmp{\mfA}$ is odd, i.e.\ not 0, 2, or 4.
  \begin{itemize}
  \item If $\prmp{\mfA}=0$, then there would by \prpref{linkH0}---part
    (f) if $\prmq{\mfA} = 1$ and part (d) if $\prmq{\mfA} = 3$---exist
    a grade $3$ perfect ideal $\mfB$ in $Q$ of class $\mathbf{H}$ with
    \begin{equation*}
      \prmm{\mfB} = \prmn{\mfA}+3\,,\quad \prmn{\mfB}=2\,,
      \qand \prmp{\mfB} \ge 1\:.
    \end{equation*}
    As $\prmm{\mfB}$ is odd, it follows from \prpref{n2} that no such
    ideal exists.
  \item If $\prmp{\mfA}=2$, then there would by \prpref{linkH2}---part
    (f) if $\prmq{\mfA} = 1$ and part (d) if $\prmq{\mfA} = 3$---exist
    a grade $3$ perfect ideal $\mfB$ in $Q$ of class $\mathbf{H}$ with
    \begin{equation*}
      \prmm{\mfB} = \prmn{\mfA}+1\,,\quad \prmn{\mfB}=2\,, 
      \qand \prmp{\mfB} \ge 1\:.
    \end{equation*}
    As $\prmm{\mfB}$ is odd, it follows from \prpref{n2} that no such
    ideal exists.

  \item If $\prmp{\mfA}=4$, then there would by \prpref{linkH2}(b)
    exist a grade $3$ perfect ideal $\mfB$ in $Q$ of class
    $\mathbf{H}$ with
    \begin{equation*}
      \prmm{\mfB} = \prmn{\mfA}+1\,,\quad \prmn{\mfB}=2\,, 
      \qand \prmp{\mfB} \ge 1\:.
    \end{equation*}
    As $\prmm{\mfB}$ is odd, it follows from \prpref{n2} that no such
    ideal exists.
  \end{itemize}

  Assume now that $\prmn{\mfA}$ is odd. The goal is to show that
  $\prmp{\mfA}$ is even, i.e.\ not 1 or 3.
  \begin{itemize}
  \item If $\prmp{\mfA}=1$, then there would by \prpref{linkH1}---part
    (f) if $\prmq{\mfA} = 1$ and part (d) if $\prmq{\mfA} = 3$---exist
    a grade $3$ perfect ideal $\mfB$ in $Q$ of class $\mathbf{H}$ with
    \begin{equation*}
      \prmm{\mfB} = \prmn{\mfA}+2\,,\quad \prmn{\mfB}=2\,, 
      \qand \prmp{\mfB} \ge 1\:.
    \end{equation*}
    As $\prmm{\mfB}$ is odd, it follows from \prpref{n2} that no such
    ideal exists.

  \item If $\prmp{\mfA}=3$, then there would by \prpref{linkH1}(b)
    exist a grade $3$ perfect ideal $\mfB$ in $Q$ of class
    $\mathbf{H}$ with
    \begin{equation*}
      \prmm{\mfB} = \prmn{\mfA}+2\,,\quad \prmn{\mfB}=2\,,
      \qand \prmp{\mfB} \ge 1\:.
    \end{equation*}
    As $\prmm{\mfB}$ is odd, it follows from \prpref{n2} that no such
    ideal exists.\qedhere
  \end{itemize}
\end{prf*}

\begin{lem}
  \label{lem:n3}
  For every grade $3$ perfect ideal $\mfB \subseteq Q$ of class
  $\mathbf{H}$ with $\prmn{\mfB} = 3$ and $\prmp{\mfB} \ge 2$ one has
  $\prmq{\mfB} \equiv_2 \prmm{\mfB}-4$.
\end{lem}

\begin{prf*}
  By \eqref{mn} and \eqref{pmqn} one has $\prmm{\mfB} \ge 4$ and
  $\prmq{\mfB} \le 3$. The assumption $\prmp{\mfB} \ge 2$ together
  with \eqref{pnqm} and \prpref{sanchez} yields $\prmp{\mfB} = 2$ or
  $\prmp{\mfB} = 4$. We treat odd and even values of $\prmm{\mfB}$
  separately.

  Assume that $\prmm{\mfB}$ is even. The goal is to show that
  $\prmq{\mfB}$ is even, i.e.\ not 1 or 3.
  \begin{itemize}
  \item If $\prmq{\mfB} = 1$, then there would by
    \prpref{linkH2}---part (f) if $\prmp{\mfB} = 2$ and part (b) if
    $\prmp{\mfB} = 4$---exist a grade $3$ perfect ideal $\mfA$ in $Q$
    of class $\mathbf{H}$ with $\prmm{\mfA} = 4$ and
    $\prmn{\mfA} = \prmm{\mfB}-3$.  As $\prmn{\mfA}$ is odd, this
    contradicts \thmref{m4}.

  \item If $\prmq{\mfB} = 3$, then there would by \prpref{linkH2}(d)
    exist a grade $3$ perfect ideal $\mfA$ in $Q$ of class
    $\mathbf{H}$ with $\prmm{\mfA} = 4$ and
    $\prmn{\mfA} = \prmm{\mfB}-3$.  As $\prmn{\mfA}$ is odd, this
    contradicts \thmref{m4}.
  \end{itemize}

  Assume now that $\prmm{\mfB}$ is odd. The goal is to show that
  $\prmq{\mfB}$ is odd, i.e.\ not 0 or 2.
  \begin{itemize}
  \item If $\prmq{\mfB} = 0$, then there would by \prpref{linkH1}(a)
    exist a grade $3$ perfect ideal $\mfA$ in $Q$ with
    \begin{equation*}
      \prmm{\mfA} = 5\,,\quad \prmn{\mfA}=\prmm{\mfB}-3\,,\quad 
      \prmp{\mfA} = 0\,, \qand \prmq{\mfA} \ge 1\:.
    \end{equation*}
    It follows from \eqref{pqr}, \pgref{ci}, and \thmref{m5} that
    $\mfA$ is of class $\mathbf{H}$. As $\prmn{\mfA}$ is even, it
    follows from \lemref{m5} that no such ideal $\mfA$ exists.

  \item If $\prmq{\mfB} = 2$, then there would by \prpref{linkH1}(a,c)
    exist a grade $3$ perfect ideal $\mfA$ in $Q$ with
    \begin{equation*}
      \prmm{\mfA} = 5\,,\quad \prmn{\mfA}=\prmm{\mfB}-3\,,\quad 
      \prmp{\mfA} = 2\,, \qand \prmq{\mfA} \ge 1\:.
    \end{equation*}
    It follows from \eqref{pqr} that $\mfA$ is of class $\mathbf{H}$.
    As $\prmn{\mfA}$ is even, it follows from \lemref{m5} that no such
    ideal $\mfA$ exists.\qedhere
  \end{itemize}
\end{prf*}

\begin{lem}
  \label{lem:n4}
  For every grade $3$ perfect ideal $\mfB \subseteq Q$ of class
  $\mathbf{H}$ with $\prmn{\mfB} = 4$ and $\prmp{\mfB} \ge 3$ one has
  $\prmq{\mfB} \equiv_2 \prmm{\mfB}-4$.
\end{lem}

\begin{prf*}
  By \eqref{mn} one has $\prmm{\mfB} \ge 4$. Notice that if
  $\prmm{\mfB} = 4$ holds, then \thmref{m4} yields
  $\prmq{\mfB} = 0 = \prmm{\mfB}-4$.  If $\prmm{\mfB} \ge 5$, then
  \prpref{linkH2}(a) yields a grade $3$ perfect ideal $\mfA$ in $Q$
  with
  \begin{equation*}
    \prmm{\mfA} = 5\,, \quad \prmn{\mfA} = \prmm{\mfB} - 3 \ge 2\,,\quad 
    \prmp{\mfA} = \prmq{\mfB}\,,\qand \prmq{\mfA} \ge 1\:.
  \end{equation*}
  It follows from \thmref{m5} and \eqref{pqr} that $\mfA$ is of class
  $\clB$ or $\mathbf{H}$.

  If $\mfB$ is of class $\clB$, then \thmref[]{m5}(a) yields
  $\prmn{\mfA} =2$, and hence $\prmm{\mfB} = 5$. Per \eqref{pqr} one
  thus has $\prmq{\mfB} = \prmp{\mfA} = 1 = \prmm{\mfB}-4$.

  If $\mfB$ is of class $\mathbf{H}$, then \lemref{m5} yields
  $\prmq{\mfB} = \prmp{\mfA} \equiv_2 \prmn{\mfA}- 1 = \prmm{\mfB}-4$.
\end{prf*}

\begin{lem}
  \label{lem:NN}
  Let $n \ge 3$ be an integer. Assume that for every grade $3$ perfect
  ideal $\mfB \subseteq Q$ of class $\mathbf{H}$ with
  $\prmn{\mfB} = n$ and $\prmp{\mfB} \ge \prmn{\mfB}-1$, one has
  $\prmq{\mfB} \equiv_2 \prmm{\mfB} - 4$. \linebreak For every grade
  $3$ perfect ideal $\mfA \subseteq Q$ of class $\mathbf{H}$ with
  $\prmm{\mfA} = n+3$ and $\prmq{\mfA} \ge \prmm{\mfA}-4$, one then
  has $\prmp{\mfA} \equiv_2 \prmn{\mfA} -1$.
\end{lem}

\begin{prf*}
  Let $n \ge 3$ be an integer and $\mfA$ a grade $3$ perfect ideal in
  $Q$ of class $\mathbf{H}$ with
  \begin{equation*}
    \tag{1}
    \prmm{\mfA} = n + 3 \ge 6 \qand \prmq{\mfA} \ge 
    \prmm{\mfA} - 4 = n - 1 \ge 2\:.
  \end{equation*}
  It needs to be shown $\prmp{\mfA}$ and $\prmn{\mfA}$ have opposite
  parity. By \eqref{pmqn} one has $\prmp{\mfA} \le n+2$; we treat the
  cases $\prmp{\mfA} \le n$, $\prmp{\mfA} = n+1$, and
  $\prmp{\mfA} = n+2$ separately.

  \emph{Case 1.} Assume that
  $0 \le \prmp{\mfA} \le n = \prmm{\mfA} - 3$ holds. As
  $\prmq{\mfA} \ge 2$ by $(1)$, there exists by \prpref{linkH0}(c) and
  $(1)$ a grade $3$ perfect ideal $\mfB$ in $Q$ with
  \begin{equation*}
    \prmm{\mfB} = \prmn{\mfA} + 3\,,\quad  \prmn{\mfB} = n\,,\quad  
    \prmp{\mfB} \ge n-1\,, \qand  \prmq{\mfB} = \prmp{\mfA}\:.  
  \end{equation*}
  For $n \ge 4$ one has $\prmq{\mfA} \ge 3$ by $(1)$, so it follows
  from \prpref{linkH0}(d) that $\mfB$ is of class $\mathbf{H}$; for
  $n=3$ the same conclusion follows from \eqref{pqr} and
  \prpref{sanchez}. The assumptions now yield the congruence
  \begin{equation*}
    \prmp{\mfA} = \prmq{\mfB} \equiv_2 \prmm{\mfB} -4 = \prmn{\mfA}-1\:.
  \end{equation*}

  \emph{Case 2.} Assume that $\prmp{\mfA} = n+1 = \prmm{\mfA} - 2$
  holds.  As $\prmq{\mfA} \ge 2$ by $(1)$, there exists by
  \prpref{linkH1}(c) and $(1)$ a grade $3$ perfect ideal $\mfB$ in $Q$
  with
  \begin{equation*}
    \prmm{\mfB} = \prmn{\mfA} + 2\,,\quad  \prmn{\mfB} = n\,,\quad  
    \prmp{\mfB} \ge n-1\,, \qand  \prmq{\mfB} = \prmp{\mfA}-1 = n\:.  
  \end{equation*}
  As $\prmq{\mfB} \ge 3$, the ideal $\mfB$ is by \eqref{pqr} of class
  $\mathbf{H}$, so the assumptions yield
  \begin{equation*}
    \prmp{\mfA} = \prmq{\mfB} + 1 \equiv_2 \prmm{\mfB} - 3 = \prmn{\mfA}-1\:.
  \end{equation*}

  \emph{Case 3.} Assume that $\prmp{\mfA} = n+2 = \prmm{\mfA}-1$
  holds.  As $\prmq{\mfA} \ge 2$ by $(1)$, there exists by
  \prpref{linkH2}(c) and $(1)$ a grade $3$ perfect ideal $\mfB$ in $Q$
  with
  \begin{equation*}
    \prmm{\mfB} = \prmn{\mfA} + 1\,,\quad  \prmn{\mfB} = n\,,\quad  
    \prmp{\mfB} \ge n-1\,, \qand  \prmq{\mfB} = \prmp{\mfA}-2=n\:.  
  \end{equation*}
  As $\prmq{\mfB} \ge 3$, the ideal $\mfB$ is by \eqref{pqr} of class
  $\mathbf{H}$, so the assumptions yield
  \begin{equation*}
    \prmp{\mfA} = \prmq{\mfB} + 2 \equiv_2 
    \prmm{\mfB} - 2 = \prmn{\mfA}-1\:.\qedhere
  \end{equation*}
\end{prf*}

\begin{lem}
  \label{lem:MM}
  Let $m \ge 7$ be an integer. Assume that for every grade $3$ perfect
  ideal $\mfA \subseteq Q$ of class $\mathbf{H}$ with
  $m-1 \le \prmm{\mfA} \le m$ and $\prmq{\mfA} \ge \prmm{\mfA}-4$ one
  has $\prmp{\mfA} \equiv_2 \prmn{\mfA} - 1$.\linebreak For every
  grade $3$ perfect ideal $\mfB \subseteq Q$ of class $\mathbf{H}$
  with $\prmn{\mfB} = m-2$ and $\prmp{\mfB} \ge \prmn{\mfB}-1$, one
  then has $\prmq{\mfB} \equiv_2 \prmm{\mfB} -4$.
\end{lem}

\begin{prf*}
  Let $m \ge 7$ be an integer and $\mfB$ a grade $3$ perfect ideal in
  $Q$ of class $\mathbf{H}$ with
  \begin{equation*}
    \tag{1}
    \prmn{\mfB} = m - 2 \ge 5 \qand \prmp{\mfB} \ge 
    \prmn{\mfB} - 1 = m - 3 \ge 4\:.
  \end{equation*}
  It needs to be shown $\prmq{\mfB}$ and $\prmm{\mfB}$ have the same
  parity.  By \eqref{pmqn} one has $\prmp{\mfB} \le \prmm{\mfB}-1$. We
  treat inequality and equality separately.

  \emph{Case 1.} Assume that $\prmp{\mfB} \le \prmm{\mfB} -2$ holds.
  As $\prmp{\mfB} \ge 4$ by $(1)$, there exists by \prpref{linkH1}(b)
  a grade $3$ perfect ideal $\mfA$ in $Q$ of class \textbf{H} with
  \begin{equation*}
    \prmm{\mfA} = m\,,\quad  \prmn{\mfA} = \prmm{\mfB}-3\,,\quad  
    \prmp{\mfA} = \prmq{\mfB}\,, \qand  \prmq{\mfA} \ge \prmp{\mfB}-1 \ge m-4\:.  
  \end{equation*}
  As $\prmq{\mfA} \ge \prmm{\mfA}-4$, the assumptions now yield
  \begin{equation*}
    \prmq{\mfB} = \prmp{\mfA} \equiv_2 \prmn{\mfA}-1 = \prmm{\mfB}-4\:.
  \end{equation*}

  \emph{Case 2.} Assume that $\prmp{\mfB} = \prmm{\mfB} - 1$ holds.
  As $\prmp{\mfB} \ge 4$ by $(1)$, there exists by \prpref{linkH2}(b)
  a grade $3$ perfect ideal $\mfA$ in $Q$ of class $\mathbf{H}$ with
  \begin{equation*}
    \prmm{\mfA} = m-1\,,\quad  \prmn{\mfA} = \prmm{\mfB}-3\,,\quad  
    \prmp{\mfA} = \prmq{\mfB}\,, \qand  \prmq{\mfA} \ge\prmp{\mfB}-2 
    \ge m-5\:.  
  \end{equation*}
  As $\prmq{\mfA} \ge \prmm{\mfA}-4$, the assumptions now yield
  \begin{equation*}
    \prmq{\mfB} = \prmp{\mfA} \equiv_2 \prmn{\mfA}-1 
    = \prmm{\mfB}-4\:. \qedhere
  \end{equation*}
\end{prf*}

\begin{prf*}[Proof of \ref{thm:pnqm-1}]
  Per \eqref{mn} one has $\prmm{\mfI} \ge 4$ and $\prmn{\mfI} \ge 2$;
  to prove the assertion we show by induction on $l \in\NN$ that the
  following hold:
  \begin{rqm}
  \item For every grade $3$ perfect ideal $\mfA \subseteq Q$ of class
    $\mathbf{H}$ with $\prmm{\mfA} = l+3$ and
    $\prmq{\mfA} \ge \prmm{\mfA} - 4$ one has
    $\prmp{\mfA} \equiv_2 \prmn{\mfA}-1$.
  \item For every grade $3$ perfect ideal $\mfB \subseteq Q$ of class
    $\mathbf{H}$ with $\prmn{\mfB} = l+1$ and
    $\prmp{\mfB} \ge \prmn{\mfB} - 1$ one has
    $\prmq{\mfB} \equiv_2 \prmm{\mfB}-4$.
  \end{rqm}

  For $l = 1$ the assertion (1) holds by \thmref{m4} and (2) holds by
  \prpref{n2}.

  For $l = 2$ the assertion (1) is \lemref{m5} and (2) is \lemref{n3}.

  Let $l = 3$. As (2) holds for $l=2$, the premise in \lemref{NN} is
  satisfied for $n=3$, and it follows that (1) holds. The assertion
  (2) is \lemref{n4}.

  Let $l \ge 4$ and assume that (1) and (2) hold for all lower values
  of $l$. By part (2) for $l-1$ the premise in \lemref{NN} is
  satisfied for $n=l$, and it follows that (1) holds for $l$. By the
  assertion (1) for $l$ and $l-1$ the premise in \lemref{MM} is
  satisfied for $m=l+3$, and it follows that (2) holds for $l$.
\end{prf*}

\section{The realizability question}
\label{sec:BGT}

\noindent
We sum up the contributions of the previous sections towards an answer
to the realizability question~\cite[Question 3.8]{LLA12}. The focus is
still on restricting the range of potentially realizable classes.

\begin{bfhpg}[Class B]
  \label{B}
  Let $\mfI \subseteq Q$ be a grade $3$ perfect ideal of class
  $\mathbf{B}$. \thmref{m4}, \pgref{ci}, and \pgref{gor} yield
  $\prmm{\mfI} \ge 5$ and $\prmn{\mfI} \ge 2$. Moreover, if
  $\prmm{\mfI}=5$ holds, then \thmref{m5} yields $\prmn{\mfI}=2$, and
  if $\prmn{\mfI}=2$ holds, then $\prmm{\mfI}$ is odd by \prpref{n2}.
\end{bfhpg}
  
  \begin{sidewaystable}
    \centering \vspace{.6\textwidth}
    \begin{tikzpicture}[x=.8cm,y=.4cm]
      \draw[white] (0,0) -- (0,-1.5); \draw[step=1.0,gray,very
      thin](0,0) grid (22,25); \fill[gray!20!white] (3,2) rectangle
      (4,3); \fill[gray!20!white] (3,7) rectangle (4,8);
      \fill[gray!20!white] (3,18) rectangle (4,19);

      \fill[gray!20!white] (4,0) rectangle (5,1); \fill[gray!20!white]
      (4,3) rectangle (6,4); \fill[gray!20!white] (4,4) rectangle
      (5,5); \fill[gray!20!white] (6,4) rectangle (7,5);
      \fill[gray!20!white] (8,6) rectangle (9,7);

      \fill[gray!20!white] (4,7) rectangle (7,8); \fill[gray!20!white]
      (5,8) rectangle (6,9); \fill[gray!20!white] (7,8) rectangle
      (8,9); \fill[gray!20!white] (4,12) rectangle (8,13);
      \fill[gray!20!white] (4,13) rectangle (5,14);
      \fill[gray!20!white] (6,13) rectangle (7,14);
      \fill[gray!20!white] (8,13) rectangle (9,14);
      \fill[gray!20!white] (4,18) rectangle (9,19);
      \fill[gray!20!white] (5,19) rectangle (6,20);
      \fill[gray!20!white] (7,19) rectangle (8,20);

      \fill[gray!20!white] (9,0) rectangle (10,2);
      \fill[gray!20!white] (10,2) rectangle (11,3);
      \fill[gray!20!white] (9,3) rectangle (11,5);
      \fill[gray!20!white] (11,3) rectangle (12,4);
      \fill[gray!20!white] (9,5) rectangle (10,6);
      \fill[gray!20!white] (11,5) rectangle (12,6);

      \fill[gray!20!white] (9,7) rectangle (12,9);
      \fill[gray!20!white] (12,7) rectangle (13,8);
      \fill[gray!20!white] (10,9) rectangle (11,10);
      \fill[gray!20!white] (12,9) rectangle (13,10);
      \fill[gray!20!white] (14,11) rectangle (15,12);

      \fill[gray!20!white] (9,12) rectangle (13,14);
      \fill[gray!20!white] (13,12) rectangle (14,13);
      \fill[gray!20!white] (9,14) rectangle (10,15);
      \fill[gray!20!white] (11,14) rectangle (12,15);
      \fill[gray!20!white] (13,14) rectangle (14,15);

      \fill[gray!20!white] (9,18) rectangle (15,19);
      \fill[gray!20!white] (9,19) rectangle (14,20);
      \fill[gray!20!white] (10,20) rectangle (11,21);
      \fill[gray!20!white] (12,20) rectangle (13,21);
      \fill[gray!20!white] (14,20) rectangle (15,21);

      \fill[gray!20!white] (15,0) rectangle (16,3);

      \fill[gray!20!white] (15,3) rectangle (17,6);
      \fill[gray!20!white] (17,4) rectangle (18,5);
      \fill[gray!20!white] (15,6) rectangle (16,7);
      \fill[gray!20!white] (17,6) rectangle (18,7);

      \fill[gray!20!white] (15,7) rectangle (18,10);
      \fill[gray!20!white] (18,8) rectangle (19,9);
      \fill[gray!20!white] (16,10) rectangle (17,11);
      \fill[gray!20!white] (18,10) rectangle (19,11);

      \fill[gray!20!white] (15,12) rectangle (19,15);
      \fill[gray!20!white] (19,13) rectangle (20,14);
      \fill[gray!20!white] (15,15) rectangle (16,16);
      \fill[gray!20!white] (17,15) rectangle (18,16);
      \fill[gray!20!white] (19,15) rectangle (20,16);
      \fill[gray!20!white] (21,17) rectangle (22,18);

      \fill[gray!20!white] (15,18) rectangle (20,21);
      \fill[gray!20!white] (20,19) rectangle (21,20);
      \fill[gray!20!white] (16,21) rectangle (17,22);
      \fill[gray!20!white] (18,21) rectangle (19,22);
      \fill[gray!20!white] (20,21) rectangle (21,22);

      \draw (3.5,3.5) node {$\mathbf{T}$}; \draw (3.5,12.5) node
      {$\mathbf{T}$}; \draw (7.5,7.5) node {$\mathbf{T}$}; \draw
      (7.5,12.5) node {$\mathbf{T}$}; \draw (7.5,18.5) node
      {$\mathbf{T}$}; \draw (12.5,7.5) node {$\mathbf{T}$}; \draw
      (12.5,12.5) node {$\mathbf{T}$}; \draw (12.5,18.5) node
      {$\mathbf{T}$}; \draw (18.5,7.5) node {$\mathbf{T}$}; \draw
      (18.5,12.5) node {$\mathbf{T}$}; \draw (18.5,18.5) node
      {$\mathbf{T}$};

      \draw (5.5,1.5) node {$\mathbf{B}$};

      \draw (10.5,4.5) node {$\mathbf{B}$}; \draw (10.5,8.5) node
      {$\mathbf{B}$}; \draw (10.5,13.5) node {$\mathbf{B}$}; \draw
      (10.5,19.5) node {$\mathbf{B}$};

      \draw (16.5,1.5) node {$\mathbf{B}$}; \draw (16.5,4.5) node
      {$\mathbf{B}$}; \draw (16.5,8.5) node {$\mathbf{B}$}; \draw
      (16.5,13.5) node {$\mathbf{B}$}; \draw (16.5,19.5) node
      {$\mathbf{B}$};

      \draw[black,thick,->](0,0) -- (23,0); \draw[black,thick,->](0,0)
      -- (0,27); \draw[black] (0,3) -- (22.5,3); \draw[black] (0,7) --
      (22.5,7); \draw[black] (0,12) -- (22.5,12); \draw[black] (0,18)
      -- (22.5,18); \draw[black] (0,25) -- (22.5,25); \draw[black]
      (4,0) -- (4,26); \draw[black] (9,0) -- (9,26); \draw[black]
      (15,0) -- (15,26); \draw[black] (22,0) -- (22,26); \draw
      (-.25,26) node {$q$}; \draw (-.25,.5) node {$0$}; \draw
      (-.25,1.75) node {$\vdots$}; \draw (-.25,2.5) node {$2$}; \draw
      (-.25,3.5) node {$0$}; \draw (-.25,4.75) node {$\vdots$}; \draw
      (-.25,5.75) node {$\vdots$}; \draw (-.25,6.5) node {$3$}; \draw
      (-.25,7.5) node {$0$}; \draw (-.25,8.75) node {$\vdots$}; \draw
      (-.25,9.5) node {$2$}; \draw (-.25,10.75) node {$\vdots$}; \draw
      (-.25,11.5) node {$4$}; \draw (-.25,12.5) node {$0$}; \draw
      (-.25,13.75) node {$\vdots$}; \draw (-.25,14.5) node {$2$};
      \draw (-.25,15.75) node {$\vdots$}; \draw (-.25,16.75) node
      {$\vdots$}; \draw (-.25,17.5) node {$5$}; \draw (-.25,18.5) node
      {$0$}; \draw (-.25,19.75) node {$\vdots$}; \draw (-.25,20.5)
      node {$2$}; \draw (-.25,21.75) node {$\vdots$}; \draw
      (-.25,22.5) node {$4$}; \draw (-.25,23.75) node {$\vdots$};
      \draw (-.25,24.5) node {$6$}; \draw (2,25.5) node {$m=4$}; \draw
      (6.5,25.5) node {$m=5$}; \draw (12,25.5) node {$m=6$}; \draw
      (18.5,25.5) node {$m=7$}; \draw (22.5,-.5) node {$p$}; \draw
      (0.5,-0.5) node {$0$}; \draw (1.5,-0.5) node {$1$}; \draw
      (2.5,-0.5) node {$2$}; \draw (3.5,-0.5) node {$3$}; \draw
      (4.5,-0.5) node {$0$}; \draw (5.5,-0.5) node {$1$}; \draw
      (6.5,-0.5) node {$2$}; \draw (7.5,-0.5) node {$3$}; \draw
      (8.5,-0.5) node {$4$}; \draw (9.5,-0.5) node {$0$}; \draw
      (10.5,-0.5) node {$1$}; \draw (11.5,-0.5) node {$2$}; \draw
      (12.5,-0.5) node {$3$}; \draw (13.5,-0.5) node {$4$}; \draw
      (14.5,-0.5) node {$5$}; \draw (15.5,-0.5) node {$0$}; \draw
      (16.5,-0.5) node {$1$}; \draw (17.5,-0.5) node {$2$}; \draw
      (18.5,-0.5) node {$3$}; \draw (19.5,-0.5) node {$4$}; \draw
      (20.5,-0.5) node {$5$}; \draw (21.5,-0.5) node {$6$};

      \draw (22,3) -- node[sloped,above] {$n=2$} (22,0); \draw (22,7)
      -- node[sloped,above] {$n=3$} (22,3); \draw (22,12) --
      node[sloped,above] {$n=4$} (22,7); \draw (22,18) --
      node[sloped,above] {$n=5$} (22,12); \draw (22,25) --
      node[sloped,above] {$n=6$} (22,18);
    \end{tikzpicture}
    \caption{Codimension $3$ Cohen--Macaulay local rings that are not
      Gorenstein have first derivation $m \ge 4$ and type $n \ge
      2$. For such rings with $m \le 7$ and $n\le 6$ the table shows
      which of the classes $\clB$, $\mathbf{H}$, and $\clT$ may exist
      according to \thmref{A}, \pgref{B}, and \pgref{T}. To avoid
      overloading the table, the possible classes $\clH{p,q}$ are only
      indicated by shading the corresponding fields in the
      $pq$-grid.}  \label{tab:2}
  \end{sidewaystable}

  \begin{bfhpg}[Class T]
    \label{T}
    Let $\mfI \subseteq Q$ be a grade $3$ perfect ideal of class
    $\mathbf{T}$; by \pgref{ci} one has $\prmm{\mfI} \ge 4$.  If
    $\prmm{\mfI}=4$ holds, then it follows from \thmref{m4} that
    $\prmn{\mfI}$ is odd and at least $3$. If $\prmm{\mfI}\ge 5$
    holds, then \prpref{n2} and \prpref{sanchez} yield
    $\prmn{\mfI} \ge 4$.
  \end{bfhpg}

\begin{bfhpg}[Summary]
  \tabref{2} illustrates which Cohen--Macaulay local rings of class
  $\mathbf{H}$ are possible per \thmref{A}. The table also shows the
  restrictions imposed by \pgref{B} and \pgref{T} on the existence of
  rings of class $\clB$ or $\clT$.
\end{bfhpg}

\begin{bfhpg}[Conjectures]
  The statement of \thmref{A} was informed by experiments conducted
  using the \textsc{Macaulay 2} implementation of Christensen and
  Veliche's classification algorithm~\cite{LWCOVl14a,M2}. Based on
  these experiments we conjecture that for $m\ge 5$ and $n \ge 3$
  there are no further restrictions on realizability of codimension
  $3$ Cohen--Macaulay local rings of class $\clB$, $\mathbf{H}$, or
  $\clT$ than those captured by \thmref[]{A}, \pgref{B}, and \pgref{T}
  and illustrated by \tabref{2}.

  Absent from \tabref{2} are rings of classes $\clC{3}$ and
  $\mathbf{G}$. A Cohen--Macaulay local ring of codimension $3$, first
  derivation $m$, and type $n$ is Gorenstein if and only if $n=1$
  holds, and they exist for odd $m \ge 3$; see \pgref{gor} and
  \pgref{reg}. For $m=3$ they are complete intersections, precisely
  the rings of class $\clC{3}$, and for $m \ge 5$ they are of class
  $\clG{m}$. Per \thmcite[3.1]{LLA12} one has $2 \le r \le m-2$ for
  rings of class $\clG{r}$ that are not Gorenstein.  Our experiments
  suggests that there are further restrictions.

  Let $\mfI \subseteq Q$ be a grade $3$ perfect ideal of class
  $\mathbf{G}$ and not Gorenstein.  By \pgref{ci}, \thmref{m4}, and
  \thmref{m5} one has $\prmm{\mfI} \ge 6$. We conjecture that the
  following hold:
  \begin{prt}
  \item If $\prmn{\mfI} = 2$, then one has
    $2 \le \prmr{\mfI} \le \prmm{\mfI}-5$ or
    $\prmr{\mfI} = \prmm{\mfI}-3$.
  \item If $\prmn{\mfI} \ge 3$, then one has
    $2 \le \prmr{\mfI} \le \prmm{\mfI}-4$.
  \end{prt}
  The next proposition proves part (b) of this conjecture for six
  generated ideals.
\end{bfhpg}

\begin{prp}
  Let $\mfA \subseteq Q$ be a grade 3 perfect ideal with
  $\prmm{\mfA}=6$ and $\prmn{\mfA} \ge 3$. If $\mfA$ is of class
  $\mathbf{G}$, then $\prmr{\mfA}=2$.
\end{prp}

\begin{prf*}
  If $\mfA$ is of class $\mathbf{G}$, then there exists by
  \prpref{linkBGT}(b) a grade $3$ perfect ideal $\mfB$ in $Q$ with
  \begin{equation*}
    \prmm{\mfB} \deq \prmn{\mfA} + 3\,, \quad \prmn{\mfB} \deq
    3\,, \qand \prmp{\mfB} \dge \min\set{\prmr{\mfA},3} \ge 2\:.
  \end{equation*}
  As $\prmm{\mfB} - \prmn{\mfB} = \prmn{\mfA} \ge 3$, \corref{MT1m2}
  yields $\prmp{\mfB} \le 2$, which forces $\prmr{\mfA}=2$.
\end{prf*}

\appendix
\section*{Appendix. Multiplication in Tor algebras of linked
  ideals---proofs}

\stepcounter{section}
\begin{stp}
  \label{stp:link}
  Let $Q$ be a local ring with maximal ideal $\mfM$.  Let
  $\mfA \subseteq Q$ be a grade $3$ perfect ideal; set
  \begin{equation*}
    m = \prmm{\mfA}\,,\quad n = \prmn{\mfA}\,,\qand l=m+n-1\:.
  \end{equation*}
  Let $A_\sbt \to Q/\mfA$ be a minimal free resolution over $Q$ and
  set
  \begin{equation*}
    \A_\sbt \deq \Tor[Q]{\sbt}{Q/\mfA}{k} \deq \H[\sbt]{\tp[Q]{A_\sbt}{k}}\:.
  \end{equation*}
  Let $\set{e_i}_{1\le i \le m}$, $\set{f_j}_{1\le j \le l}$, and
  $\set{g_h}_{1\le h \le n}$ denote bases for $A_1$, $A_2$, and $A_3$.
\end{stp}

\begin{bfhpg}[Multiplicative structures]
  Adopt \stpref{link}. The homology classes of the bases for $A_1$,
  $A_2$, and $A_3$ yield bases
  \begin{equation}
    \label{eq:Abasis}
    \e_1,\ldots,\e_m \ \text{ for $\A_1$}\,,\quad
    \f_1,\ldots, \f_{l} \ \text{ for $\A_2$}\,,\qand
    \g_1,\ldots,\g_n \ \text{ for $\A_3$}\:.
  \end{equation}
  As recalled in \pgref{ms}, bases can be chosen such that the
  multiplication on $\A_\sbt$ is one of:
  \begin{equation}
    \label{eq:Aefg}
    \begin{aligned}
      \textbf{C}(3): \quad & \e_1\e_2 = \f_3 \ \ \e_2\e_3 = \f_1 \ \
      \e_3\e_1 = \f_2
      & \ \e_i\f_i = \g_1 \ \text{ for } \ 1\le i \le 3\\[.5ex]
      \textbf{T}: \quad & \e_1\e_2 = \f_3 \ \ \e_2\e_3 = \f_1 \ \
      \e_3\e_1  = \f_2\\[.5ex]
      \textbf{B}: \quad & \e_1\e_2 = \f_3 \
      & \ \e_i\f_i = \g_1 \ \text{ for } \ 1\le i \le 2\\[.5ex]
      \clG{r}: \quad & [r\ge 2]
      & \e_i\f_i = \g_1 \ \text{ for } \ 1\le i \le r\\[.5ex]
      \clH{p,q}: \quad & \e_{p+1}\e_i = \f_i \ \text{ for } \ 1\le
      i\le p & \e_{p+1}\f_{p+j} = \g_j \ \text{ for } \ 1\le j\le q
    \end{aligned}
  \end{equation}
\end{bfhpg}
\begin{bfhpg}[Linkage]
  \label{link}
  Adopt \stpref{link}.  Let $\mfX \subset \mfA$ be a complete
  intersection ideal generated by a regular sequence $x_1, x_2, x_3$
  and set $\mfB = (\mfX:\mfA)$. Let $K_\sbt = \Kzl{x_1, x_2, x_3}$ be
  the Koszul complex and $\mapdef{\grf_\sbt}{K_\sbt}{A_\sbt}$ be a
  lift of the canonical surjection $Q/\mfX \to Q/\mfA$ to a morphism
  of DG algebras.
  \begin{equation*}
    \begin{gathered}
      \label{diag*}
      \xymatrix@=2pc{ K_\sbt \ar@{->}[d]^{\grf_\sbt} & 0\ar@{->}[r] &
        K_3\ar@{->}[r]^{\dif[3]{K_\sbt}}\ar@{->}[d]^{\grf_3}
        &K_2\ar@{->}[r]^{\dif[2]{K_\sbt}}\ar@{->}[d]^{\grf_2}
        &K_1\ar@{->}[r]^{\dif[1]{K_\sbt}}\ar@{->}[d]^{\grf_1}
        &K_0\ar@{->}[r]\ar@{=}[d]^{1_Q} &0 \\
        A_\sbt & 0\ar@{->}[r] &A_3\ar@{->}[r]^{\dif[3]{A_\sbt}}
        &A_2\ar@{->}[r]^{\dif[2]{A_\sbt}} &
        A_1\ar@{->}[r]^{\dif[1]{A_\sbt}} &A_0\ar@{->}[r] &0 }
    \end{gathered}
  \end{equation*}
  Let $\ee_1,\ee_2,\ee_3$ denote the generators of $K_1$. From the
  mapping cone of this morphism one gets, see e.g.\
  \prpcite[1.6]{AKM-88}, a free resolution of $Q/\mfB$ over $Q$:
  \begin{equation*}
    D_\sbt \deq 0 \lra A_1^* \lra A_2^* \oplus K_2 \lra A_3^*\oplus K_1 \lra Q\;.
  \end{equation*}
  The complex $D_\sbt$ carries a multiplicative structure given by
  \thmcite[1.13]{AKM-88} in terms of the basis
  \begin{gather}
    \label{eq:E}
    \begin{pmatrix}g_1^*\\0\end{pmatrix},\ldots,
    \begin{pmatrix}g_n^*\\0\end{pmatrix},\;
    \begin{pmatrix}0\\\ee_1\end{pmatrix},\;
    \begin{pmatrix}0\\\ee_2\end{pmatrix},\;
    \begin{pmatrix}0\\\ee_3\end{pmatrix}
    \quad\text{for $D_1\:,$}\\[1ex]
    \label{eq:F}
    \begin{pmatrix}f_1^*\\0\end{pmatrix},\ldots,
    \begin{pmatrix}f_{l}^*\\0\end{pmatrix},\;
    \begin{pmatrix}0\\\ee_2\ee_3\end{pmatrix},\;
    \begin{pmatrix}0\\\ee_1\ee_3\end{pmatrix},\;
    \begin{pmatrix}0\\\ee_1\ee_2\end{pmatrix}
    \quad\text{for $D_2\;,$ and}\\[1ex]
    \label{eq:G}
    e_1^*,\ldots, e^*_m \quad\text{for $D_3\:.$}
  \end{gather}
  Here we recall from \thmcite[1.13]{AKM-88} the products that involve
  the last three vectors in the basis for $D_1$.
  \begin{gather}
    \label{eq:10D}
    \binom{0}{\ee_i}\binom{0}{\ee_j} \deq
    -\binom{0}{\ee_i\ee_j} \quad\text{for } i,j\in\set{1,2,3}\:.\\[1ex]
    \label{eq:11D}
    \binom{0}{\ee_i}\binom{g_j^*}{0} \deq
    \sum_{h=1}^{l} g_j^*(\grf_1(\ee_i)f_h)\binom{f_h^*}{0}\:.\\[1ex]
    \label{eq:21D}
    \begin{aligned}
      \binom{0}{\ee_i}\binom{0}{\ee_i\ee_j} &\deq 0 \quad\text{and}\\
      \binom{0}{\ee_h}\binom{0}{\ee_i\ee_j} &\deq \pm\dif[1]{A_\sbt}
      \in A_1^* \quad\text{for } \set{h,i,j} = \set{1,2,3}\,.
    \end{aligned}\\
    \label{eq:12D}
    \binom{0}{\ee_i}\binom{f_j^*}{0} \deq \sum_{h=1}^n
    f_j^*(\grf_1(\ee_i)e_h)e_h^*\:.
  \end{gather}

  From the free resolution $D_\sbt$ one gets a minimal free resolution
  of $Q/\mfB$ that we denote $B_\sbt$.  From \cite[(1.80)]{AKM-88} one
  has :
  \begin{align*}
    \rnk[Q] B_1 
    &= n+3-\rnk[k]{(\grf_3\otimes k)}-\rnk[k]{(\grf_2\otimes k)}\:,\\
    \rnk[Q] B_2 
    &= l+3-\rnk[k]{(\grf_2\otimes k)}-\rnk[k]{(\grf_1\otimes k)}\:,\text{ and}\\
    \rnk[Q] B_3 &= m -\rnk[k]{(\grf_1\otimes k)}\:.
  \end{align*}
  The identity $1 - \rnk[Q]{B_1} + \rnk[Q]{B_2} - \rnk[Q]{B_3} = 0$
  yields $\rnk[k]{(\grf_3\otimes k)}=0\:.$

  Set
  \begin{equation*}
    \B_\sbt \deq \Tor[Q]{\sbt}{Q/\mfB}{k} \deq \H[\sbt]{\tp[Q]{B_\sbt}{k}} 
    \deq \H[\sbt]{\tp[Q]{D_\sbt}{k}}
  \end{equation*}
  and denote the homology classes of the basis vectors as follows:
  \begin{align*}
    \E_1,\ldots, \E_n, \E_{n+1},\E_{n+2},\E_{n+3} 
    &\quad\text{for the vectors in \eqref{E},}\\
    \F_1,\ldots, \F_l, \F_{l+1},\F_{l+2},\F_{l+3} 
    &\quad\text{for the vectors in \eqref{F}, and}\\
    \G_1,\ldots, \G_m &\quad\text{for the vectors in \eqref{G}.}
  \end{align*}
  Notice that \eqref{21D} yields
  \begin{equation}
    \label{eq:zero}
    \E_{n+i}\F_{l+j} \deq 0 \quad\text{for } i,j\in\set{1,2,3}\:.
  \end{equation}
\end{bfhpg}

\begin{bfhpg}[Linkage via minimal generators]
  \label{link1}
  Adopt the setup established in \stpref[]{link}--\pgref{link}.  If
  $x_1,x_2,x_3$ are part of a minimal system of generators for $\mfA$,
  then the homomorphism $\mapdef{\grf_1}{K_1}{A_1}$ is given by
  $\grf_1(\ee_i) = e_i$ for $1\leq i\leq 3$.
  
  One has $\rnk[k](\grf_1\otimes k)=3$. More precisely, let $h$, $i$,
  and $j$ denote the three elements in $\set{1,2,3}$, now
  \prpcite[1.6]{AKM-88} yields
  \begin{equation*}
    \dif[3]{D_\sbt}(e_i^*) \deq \pm\binom{0}{\ee_h\ee_j} \mod \mfM D_2\:,
  \end{equation*}
  so in the algebra $\B_\sbt$ one has
  \begin{equation}
    \label{eq:rank}
    \F_{l+1}=\F_{l+2}=\F_{l+3} = 0 \qqand \G_1 = \G_2 = \G_3 = 0\:.
  \end{equation}
  Moreover, the rank of $\grf_2\otimes k$ is the number of linearly
  independent products $e_ie_j$ mod $\mfM A_2$ for
  $i,j\in\set{1,2,3}$. More precisely, if
  $e_ie_j = \pm f_k \mod \mfM A_2$ holds for some
  $k\in\set{1,\ldots,l}$, then \prpcite[1.6]{AKM-88} yields
  \begin{equation*}
    \dif[2]{D_\sbt}\binom{f_k^*}{0} \deq \pm\binom{0}{\ee_h} \mod \mfM D_1\:,
  \end{equation*}
  so in $\B_\sbt$ one has
  \begin{equation}
    \label{eq:split}
    \E_{n+h} = 0 \qand \F_k = 0\:.
  \end{equation}

  In particular, it follows from \eqref{rank} that the only possible
  non-zero products among the basis vectors for $\B_\sbt$ that involve
  $\E_{n+1}$, $\E_{n+2}$, or $\E_{n+3}$ are:
  \begin{align}
    \label{eq:11}
    \E_{n+i}\E_{j} &\deq
                     \sum_{h=1}^{l} \g_j^*(\e_i\f_h)\F_h\quad\text{for } \ 
                     1 \le i \le 3 \ \text{ and } \ 1 \le j \le n\:.\\[1ex]
    \label{eq:12}
    \E_{n+i}\F_j &\deq \sum_{h=4}^n
                   \f_j^*(\e_i\e_h)\G_h\quad\text{for } \ 
                   1 \le i \le 3 \ \text{ and } \ 1 \le j \le l\:.
  \end{align}
\end{bfhpg}

\begin{bfhpg}[Regular sequences]
  Adopt \stpref{link}. Let $x_1,\ldots,x_m$ be a minimal set of
  generators for $\mfA$; the homomorphism
  $\mapdef{\dif[1]{A_\sbt}}{A_1}{A_0=Q}$ is given by
  $\dif[1]{A}(e_i)=x_i$.  By a standard argument, see for example the
  proofs of \prpcite[2.2]{CVW-2} or \lemcite[8.2]{DABDEs74}, one can
  add elements from $\mfM\mfA$ to modify the sequence of generators
  such that $x_1, x_2, x_3$ form a regular sequence. This corresponds
  to a change of basis on $A_1$, where $e_i$ for $1\le i\le 3$ is
  replaced by $e_i + \sum_{j=1}^m a_{ij}e_j$ with coefficients
  \mbox{$a_{ij}\in \mfM$}. Notice that the homology classes in
  $\A_\sbt$ of the basis vectors do not change: one has
  $[e_i + \sum_{j=1}^m a_{ij}e_j] = [e_i] = \e_i$. This means that
  given any basis $\e_1,\ldots,\e_m$ for $\A_1$ one can without loss
  of generality assume that the generators $x_1,x_2,x_3$ of $\mfA$
  corresponding to $\e_1,\e_2,\e_3$ form a regular sequence. We
  tacitly do so in \pgref{prf:linkBGT}--\pgref{prf:linkH2}.

  For ease of reference, still in
  \pgref{prf:linkBGT}--\pgref{prf:linkH2}, we recall from \eqref{pqr}
  the parameters that describe the multiplicative structures:
  \begin{equation}
    \label{eq:pqrmn}
    \begin{array}{r|ccc}
      \text{Class of $\mfA$} & \prmp{\mfA} & \prmq{\mfA} & \prmr{\mfA}\\
      \hline
      \clB & 1 &1 &2 \\
      \clC{3} & 3 &1 &3 \\
      \clG{r}\ [r\ge 2]& 0 &1 &r \\
      \clH{p,q} & p & q & q\\
      \clT & 3 &0 &0 \\
    \end{array}
  \end{equation}  
\end{bfhpg}

The rest of this appendix is taken up by the arguments for
\prpref[Propositions~]{linkBGT}--\prpref[]{linkH2}.

\begin{bfhpg}[Proof of Proposition \ref{prp:linkBGT}]
  \label{prf:linkBGT}
  Adopt \stpref{link}.

  (a): One may assume that the nonzero products of elements from
  \eqref{Abasis} are
  \begin{equation*}
    \e_1\e_2 = \f_3\ \qqand \e_1\f_1 = \g_1 = \e_2\f_2\:.
  \end{equation*}
  Proceeding as in \pgref{link}, it follows from \pgref{link1} that
  $\mfA$ is directly linked to an ideal $\mfB$ whose Tor algebra
  $\B_\sbt$ has bases
  \begin{equation*}
    \begin{aligned}
      \E_1,\ldots,\E_n,\E_{n+1},\E_{n+2} \quad &\text{for } \B_1\,,\\
      \F_1,\F_2,\F_4,\ldots, \F_{l}  \quad &\text{for } \B_2\,,\text{ and}\\
      \G_4,\ldots,\G_m \quad &\text{for } \B_3\:.
    \end{aligned}
  \end{equation*}
  In particular, one has $\prmm{\mfB} = n + 2$ and
  $\prmn{\mfB} = m-3$. Further it follows from \eqref{11} and
  \eqref{12} that the nonzero products in $\B_{\sbt}$ that involve
  $\E_{n+1}$ and $\E_{n+2}$ are precisely
  \begin{equation*}
    \tag{1}\label{eq:BGHT-1}
    \E_{n+1}\E_1 = \F_{1} \qqand \E_{n+2}\E_1 = \F_{2}\:.
  \end{equation*}
  This means that $\prmp{\mfB}$ is at least $2$, so $\mfB$ is of class
  $\mathbf{H}$ or $\clT$; see \eqref{pqrmn}.

  Assume towards a contradiction that $\mfB$ is of class $\clT$. Per
  \eqref{Aefg} there are bases $\set{\E'_i}$ for $\B_1$ and
  $\set{\F'_j}$ for $\B_2$ with nonzero products
  \begin{equation*}
    \E'_{1}\E'_2 = \F'_3\,,\quad \E'_{2}\E'_3 = \F'_1\,,
    \qand \E'_{3}\E'_1 = \F'_2\:.
  \end{equation*}
  Write $\E_{n+1} = \sum_{i=1}^{n+2}\gra_{i}\E'_i$ and
  $\E_{n+2} = \sum_{i=1}^{n+2}\grb_{i}\E'_i$; now \eqref{BGHT-1}
  yields
  \begin{equation*}
    \F_1 \deq \sum_{i=1}^3\gra_i\E'_i\E_1 \qqand 
    \F_2 \deq \sum_{i=1}^3\grb_i\E'_i\E_1\:.
  \end{equation*}
  As the vectors $\F_1$ and $\F_2$ are linearly independent, so are
  the vectors $(\gra_1,\gra_2, \gra_3)$ and $(\grb_1,\grb_2,
  \grb_3)$. That is, the matrix
  \begin{equation*}
    \begin{pmatrix}
      \gra_1 & \gra_2 & \gra_3\\ \grb_1 & \grb_2 & \grb_3
    \end{pmatrix}
  \end{equation*}
  has rank $2$, whence it follows that the product
  \begin{align*}
    \E_{n+1}\E_{n+2} = (\gra_2\grb_3 - \gra_3\grb_2)\F'_1
    + (\gra_3\grb_1 - \gra_1\grb_3)\F'_2
    + (\gra_1\grb_2 - \gra_2\grb_1)\F'_3
  \end{align*}
  is nonzero, and that contradicts \eqref{BGHT-1}.

  (b): One may assume that the nonzero products of elements from
  \eqref{Abasis} are
  \begin{equation*}
    \e_i\f_i = \g_1 \quad\text{for } 1\le i \le \prmr{\mfA}\:.
  \end{equation*}
  Proceeding as in \pgref{link}, it follows from \pgref{link1} that
  $\mfA$ is directly linked to an ideal $\mfB$ whose Tor algebra
  $\B_\sbt$ has bases
  \begin{equation*}
    \tag{2}\label{eq:BGHT-2}
    \begin{aligned}
      \E_1,\ldots,\E_n,\E_{n+1},\E_{n+2},\E_{n+3} \quad &\text{for } \B_1\,,\\
      \F_1,\ldots, \F_{l}  \quad &\text{for } \B_2\,,\text{ and}\\
      \G_4,\ldots,\G_m \quad &\text{for } \B_3\:.
    \end{aligned}
  \end{equation*}
  In particular, one has $\prmm{\mfB} = n + 3$ and
  $\prmn{\mfB} = m-3$.  Further it follows from \eqref{11} and
  \eqref{12} that the nonzero products in $\B_{\sbt}$ that involve
  $\E_{n+1}$, $\E_{n+2}$, and $\E_{n+3}$ are precisely
  \begin{equation*}
    \E_{n+i}\E_1 = \F_i \quad\text{for } 1 \le i \le \min\set{\prmr{\mfA},3}\:.
  \end{equation*}
  As $\prmr{\mfA} \ge 2$ the ideal $\mfB$ is per \eqref{pqrmn} of
  class $\mathbf{H}$ or $\clT$, and the argument from the proof of
  part (a) applies to show that $\mfB$ is not of class $\clT$.

  (c): One may assume that the nonzero products of elements from
  \eqref{Abasis} are
  \begin{equation*}
    \e_3\e_4 = \f_5\,, \quad \e_4\e_5 = \f_3\,, \qand \e_5\e_3 = \f_4\:.
  \end{equation*} 
  Proceeding as in \pgref{link}, it follows from \pgref{link1} that
  $\mfA$ is directly linked to an ideal $\mfB$ whose Tor algebra
  $\B_\sbt$ has basis \eqref{BGHT-2}. In particular, one has
  $\prmm{\mfB} = n + 3$ and $\prmn{\mfB} = m-3$.  Further it follows
  from \eqref{11} and \eqref{12} that the nonzero products in
  $\B_{\sbt}$ that involve $\E_{n+1}$, $\E_{n+2}$, and $\E_{n+3}$ are
  precisely
  \begin{equation*}
    \E_{n+3}\F_5 = \G_{4} \qqand \E_{n+3}\F_4 = -\G_{5}\:. \qedhere
  \end{equation*}
  Thus one has $\prmq{\mfB} \ge 2$, in particular $\mfB$ is of class
  $\mathbf{H}$; see \eqref{pqrmn}.

  (d): One may assume that the nonzero products of elements from
  \eqref{Abasis} are
  \begin{equation*}
    \e_2\e_3 = \f_4\,, \quad \e_3\e_4 = \f_3,\, \qand \e_4\e_2 = \f_3\:.
  \end{equation*} 
  Proceeding as in \pgref{link}, it follows from \pgref{link1} that
  $\mfA$ is directly linked to an ideal $\mfB$ whose Tor algebra
  $\B_\sbt$ has bases
  \begin{equation*}
    \begin{aligned}
      \E_1,\ldots,\E_n,\E_{n+2},\E_{n+3} \quad &\text{for } \B_1\,,\\
      \F_1,\F_2,\F_3, \F_5,\ldots, \F_{l} \quad
      &\text{for } \B_2\,,\text{ and}\\
      \G_4,\ldots,\G_m \quad &\text{for } \B_3\:.
    \end{aligned}
  \end{equation*}
  In particular, one has $\prmm{\mfB} = n + 2$ and
  $\prmn{\mfB} = m-3$.  Further it follows from \eqref{11} and
  \eqref{12} that the nonzero products in $\B_{\sbt}$ that involve
  $\E_{n+2}$ and $\E_{n+3}$ are precisely
  \begin{equation*}
    \tag{3}\label{eq:BGHT-3}
    \E_{n+3}\F_3 = \G_{4} = - \E_{n+2}\F_2\:.
  \end{equation*}
  As $\prmn{\mfA} \ge 2$ by \pgref{gor}, one has $\prmm{\mfB} \ge 4$,
  $\prmq{\mfB} \ge 1$, and $\prmr{\mfB} \ge 2$, so per \eqref{pqrmn}
  the ideal $\mfB$ is of class $\clB$, $\mathbf{G}$, or $\mathbf{H}$,
  and if $\mfB$ is of class $\clB$ or $\mathbf{G}$, then
  $\prmq{\mfB}=1$ holds. If $\mfB$ is of class $\mathbf{H}$, then per
  \eqref{Aefg} there are bases $\set{\E'_i}$ for $\B_1$, $\set{\F'_j}$
  for $\B_2$, and $\set{\G'_k}$ for $\B_3$ with nonzero products
  \begin{equation*}
    \begin{alignedat}{2}
      \E'_{p'+1}\E'_j &= \F'_{j}& &\quad\text{for } 1 \le j \le p'
      \quad\text{and}\\
      \E'_{p'+1}\F'_{p'+i} &= \G'_{i}& &\quad\text{for } 1 \le i \le
      q'\:.
    \end{alignedat}
  \end{equation*}
  Write $\E_{n+3} = \sum_{i=1}^{n+2}\gra_{i}\E'_i$. It follows from
  the nonzero product $\E_{n+3}\F_3 = \G_{4}$ that $\gra_{p'+1}$ is
  nonzero. As $q' = \prmr{\mfB} \ge 2$ per \eqref{pqrmn} it follows
  that there are two linearly independent products of the form
  $\E_{n+3}\F_j$ which contradicts \eqref{BGHT-3}.

  (e): One may assume that the nonzero products of elements from
  \eqref{Abasis} are
  \begin{equation*}
    \e_1\e_2 = \f_3\,, \quad \e_2\e_3 = \f_1\,, \qand \e_3\e_1 = \f_2\:.
  \end{equation*} 
  Proceeding as in \pgref{link}, it follows from \pgref{link1} that
  $\mfA$ is directly linked to an ideal $\mfB$ whose Tor algebra
  $\B_\sbt$ has bases
  \begin{equation*}
    \begin{aligned}
      \E_1,\ldots,\E_n \quad &\text{for } \B_1\,,\\
      \F_4,\ldots, \F_{l} \quad
      &\text{for } \B_2\,,\text{ and}\\
      \G_4,\ldots,\G_m \quad &\text{for } \B_3\:.
    \end{aligned}
  \end{equation*}with
  $\prmm{\mfB} = n$ and $\prmn{\mfB} = m$.
\end{bfhpg}

\begin{bfhpg}[Proof of Proposition \ref{prp:linkH0}]
  \label{prf:linkH0}
  Adopt \stpref{link} and set $p = \prmp{\mfA}$ and $q = \prmq{\mfA}$.
  One may assume that the nonzero products of elements from
  \eqref{Abasis} are
  \begin{equation*}
    \e_1\e_{3+i} = \f_{q+i} \quad\text{for } 1 \le i \le p\qqand
    \e_1\f_j = \g_j \quad\text{for } 1 \le j \le q\:.
  \end{equation*}
  Proceeding as in \pgref{link}, it follows from \pgref{link1} that
  $\mfA$ is directly linked to an ideal $\mfB$ whose Tor algebra
  $\B_\sbt$ has bases
  \begin{equation*}
    \begin{aligned}
      \E_1,\ldots,\E_n,\E_{n+1},\E_{n+2},\E_{n+3} \quad &\text{for } \B_1\,,\\
      \F_1,\ldots, \F_{l}  \quad &\text{for } \B_2\,,\text{ and}\\
      \G_4,\ldots,\G_m \quad &\text{for } \B_3\:.
    \end{aligned}
  \end{equation*}
  In particular, one has $\prmm{\mfB} = \prmn{} + 3$ and
  $\prmn{\mfB} = \prmm{}-3$. Further it follows from \eqref{11} and
  \eqref{12} that the nonzero products in $\B_{\sbt}$ that involve
  $\E_{n+1}$ are precisely
  \begin{equation*}
    \tag{1}\label{eq:H01-1}
    \begin{alignedat}{2}
      \E_{n+1}\E_j &= \F_{j}& &\quad\text{for } 1 \le j \le q\quad\text{and}\\
      \E_{n+1}\F_{q+i} &= \G_{3+i}& &\quad\text{for } 1 \le i \le p\:.
    \end{alignedat}
  \end{equation*}
  Set $p' = \prmp{\mfB}$ and $q' = \prmq{\mfB}$; evidently one has
  \begin{equation*}
    \tag{2}\label{eq:H01-2}
    p' \ge q \qqand q' \ge p\:.
  \end{equation*}
  Notice that $\mfB$ is not complete intersection, as one has
  $\prmm{\mfB} \ge 4$; cf.~\pgref{ci}. That is, $\mfB$ is of class
  $\clB$, $\mathbf{G}$, $\mathbf{H}$, or $\clT$.

  (a): If $p \ge 1$ holds, then \eqref{H01-2} yields $q' \ge 1$ and
  per \eqref{pqrmn} it follows that the ideal $\mfB$ is of class
  $\mathbf{B}$, $\mathbf{G}$, or $\mathbf{H}$.

  If $\mfB$ is of class $\clB$, then there are per \eqref{Aefg} bases
  $\set{\E'_i}$ for $\B_1$, $\set{\F'_j}$ for $\B_2$, and
  $\set{\G'_k}$ for $\B_3$ with nonzero products
  \begin{equation*}
    \tag{3}\label{eq:H01-3}
    \E'_{1}\E'_2 = \F'_3 \qqand \E'_{1}\F'_1 = \G'_1 = \E'_{2}\F'_2\:.
  \end{equation*}
  By \eqref{pqrmn} and \eqref{H01-2} one has $1 = p' \ge q$, so assume
  towards a contraction that $q=0$ holds.  Write
  $\E_{n+1} = \sum_{i=1}^{n+3}\gra_{i}\E'_i$ and
  $\F_{q+1} = \sum_{j=1}^{l}\grb_{j}\F'_j$. By \eqref{H01-1} and
  \eqref{H01-3} one has
  \begin{equation*}
    0 \ne \G_4 = \E_{n+1}\F_{q+1} = 
    (\gra_1\grb_1 + \gra_{2}\grb_{2})\G'_1
  \end{equation*}
  and, therefore, $\gra_{1} \ne 0$ or $\gra_{2} \ne 0$. From the
  equalities
  \begin{equation*}
    \E_{n+1}\E'_{1} = \gra_{2}\E'_{2}\E'_{1} = -\gra_{2}\F'_3 \qqand 
    \E_{n+1}\E'_2  = \gra_{1}\E'_{1}\E'_2 = \gra_{1}\F'_3
  \end{equation*}
  it follows that $\E_{n+1}\E_j$ is nonzero for some $j$, which
  contradicts the assumption $q=0$; see~\eqref{H01-1}. Thus
  $q = 1 = p'$ holds.

  If $\mfB$ is of class $\mathbf{G}$, then \eqref{pqrmn} and
  \eqref{H01-2} yield $0 = p' \ge q$, so $p'=0=q$ holds.

  If $\mfB$ is of class $\mathbf{H}$, then there are per \eqref{Aefg}
  bases $\set{\E'_i}$ for $\B_1$, $\set{\F'_j}$ for $\B_2$, and
  $\set{\G'_k}$ for $\B_3$ with nonzero products
  \begin{equation*}
    \tag{4}\label{eq:H01-4}
    \begin{alignedat}{2}
      \E'_{p'+1}\E'_j &= \F'_{j}& &\quad\text{for } 1 \le j \le p'\quad\text{and}\\
      \E'_{p'+1}\F'_{p'+i} &= \G'_{i}& &\quad\text{for } 1 \le i \le
      q'\:.
    \end{alignedat}
  \end{equation*}
  Write $\E_{n+1} = \sum_{i=1}^{n+3}\gra_{i}\E'_i$. The product
  $\E_{n+1}\F_{q+1} = \G_4$ is by \eqref{H01-1} nonzero, so
  \eqref{H01-4} yields $\gra_{p'+1} \ne 0$. From
  \begin{equation*}
    \E_{n+1}\E'_j  = \gra_{p'+1}\E'_{p'+1}\E'_j = \gra_{p'+1}\F'_j \ne 0
    \quad\text{for } 1\le j\le p'
  \end{equation*}
  it follows that there are $p'$ linearly independent products of the
  form $\E_{n+1}\E_j$, which forces $p' \le q$; see
  \eqref{H01-1}. Thus $p'=q$ holds by \eqref{H01-2}.

  (b): If $p \ge 2$ holds, then \eqref{H01-2} yields $q' \ge 2$, and
  per \eqref{pqrmn} it follows that the ideal $\mfB$ is of class
  $\mathbf{H}$. Moreover, $p'=q$ holds by (a).
  
  (c): Assume that $q \ge 2$ holds. As one has $p'\ge q$, see
  \eqref{H01-2}, it follows per \eqref{pqrmn} that $\mfB$ is of class
  $\mathbf{H}$ or $\mathbf{T}$. By \eqref{H01-2} it suffices to prove
  that $q' \le p$ holds.

  If $\mfB$ is of class $\clT$, then per \eqref{pqrmn} one has $0=q'$
  so $q'\le p$ trivially holds.

  If $\mfB$ is of class $\mathbf{H}$, then as in the proof of part (a)
  there exist bases with the nonzero products given in \eqref{H01-4}.
  Write
  \begin{equation*}
    \E_{n+1} \deq \sum_{i=1}^{n+3}\gra_{i}\E'_i\:,\quad
    \E_{1} \deq \sum_{i=1}^{n+3}\grb_{i}\E'_i\:, \qand
    \E_{2} \deq \sum_{i=1}^{n+3}\grg_{i}\E'_i\:.
  \end{equation*}
  Now \eqref{H01-1} and \eqref{H01-4} yield
  \begin{align*}
    0 \ne \F_1 = \E_{n+1}\E_1 
    &\deq \sum_{j=1}^{p'}(\gra_{p'+1}\grb_{j} 
      - \gra_j\grb_{p'+1})\F'_j\quad\text{and}\\
    0 \ne \F_2 = \E_{n+1}\E_2 
    &\deq \sum_{j=1}^{p'}(\gra_{p'+1}\grg_{j} - \gra_j\grg_{p'+1})\F'_j\:.
  \end{align*}
  The vectors $\F_1$ and $\F_2$ are linearly independent, while one
  has
  \begin{equation*}
    \grg_{p'+1}\sum_{j=1}^{p'}\gra_j\grb_{p'+1}\F'_j 
    \deq \grb_{p'+1}\sum_{i=j}^{p'}\gra_j\grg_{p'+1}\F'_j\:;
  \end{equation*}
  it follows that $\gra_{p'+1}$ is nonzero.  Thus, one has
  \begin{equation*}
    \E_{n+1}\F'_{p'+i} = \gra_{p'+1}\E'_{p'+1}\F'_{p'+i} =
    \gra_{p'+1}\G'_i \quad\text{for } 1 \le i \le q'\:.
  \end{equation*}
  It follows that there are $q'$ linearly independent products of the
  form $\E_{n+1}\F_{j}$, which forces $q'\le p$; see \eqref{H01-1}.

  (d): If $q \ge 3$ holds, then \eqref{H01-2} yields $p'\ge 3$, and
  per \eqref{pqrmn} it follows that the ideal $\mfB$ is of class
  $\mathbf{H}$ or $\clT$. If $\mfB$ is of class $\clT$, then there are
  per \eqref{Aefg} bases $\set{\E'_i}$ for $\B_1$ and $\set{\F'_j}$
  for $\B_2$ with nonzero products
  \begin{equation*}
    \tag{5}\label{eq:H01-5}
    \E'_{1}\E'_2 = \F'_3\,,\quad \E'_{2}\E'_3 = \F'_1\,,
    \qand \E'_{3}\E'_1 = \F'_2\:.
  \end{equation*}
  It follow that for any element $\E'$ in $\B_1$, the map
  $\B_1 \to \B_2$ given by multiplication by $\E'$ has rank at most
  $2$. However, multiplication by $\E_{n+1}$ has rank $q \ge 3$; see
  \eqref{H01-1}. Thus $\mfB$ is not of class $\clT$ and hence of class
  $\mathbf{H}$; finally (c) yields $q'=p$.

  (e): Assume that $q = 1$ holds. By \eqref{H01-2} one has $p' \ge 1$,
  so the ideal $\mfB$ is per \eqref{pqrmn} of class $\clB$,
  $\mathbf{H}$, or $\clT$.  If $\mfB$ is of class $\clT$ then, as in
  the proof of part (d), there are bases with nonzero products as in
  \eqref{H01-5}.  Write $\E_{n+1} = \sum_{i=1}^{n+3}\gra_{i}\E'_i$ and
  $\E_{1} = \sum_{i=1}^{n+3}\grb_{i}\E'_i$. By \eqref{H01-1} and
  \eqref{H01-5} one has
  \begin{equation*}
    0 \ne \F_1 = \E_{n+1}\E_1 = (\gra_2\grb_3 - \gra_3\grb_2)\F'_1
    + (\gra_3\grb_1 - \gra_1\grb_3)\F'_2
    + (\gra_1\grb_2 - \gra_2\grb_1)\F'_3\:.
  \end{equation*}
  It follows that at least one of $\gra_1$, $\gra_2$, or $\gra_3$ is
  nonzero. From the expressions
  \begin{align*}
    \E_{n+1}\E'_1 & \deq \gra_{2}\E'_{2}\E'_1 + \gra_3\E'_3\E'_1  
                    \deq \gra_{3}\F'_2 - \gra_2\F'_3\:, \\ 
    \E_{n+1}\E'_2 & \deq \gra_{1}\E'_{1}\E'_2 + \gra_3\E'_3\E'_2  
                    \deq \gra_{1}\F'_3 - \gra_3\F'_1\:, \ \text{ and} \\ 
    \E_{n+1}\E'_3 & \deq \gra_{1}\E'_{1}\E'_3 + \gra_{2}\E'_{2}\E'_3 
                    \deq \gra_2\F'_1 -\gra_{1}\F'_2
  \end{align*}
  it now follows that there are two linearly independent products of
  the form $\E_{n+1}\E_j$, which contradicts the assumption $q=1$; see
  \eqref{H01-1}. Thus $\mfB$ is not of class $\clT$.

  (f): It follows from (e) that the ideal $\mfB$ is of class $\clB$ or
  $\mathbf{H}$.  If $\mfB$ is of class $\clB$, then as in the proof of
  part (a) there exist bases with the nonzero products given in
  \eqref{H01-3}.  Write $\E_{n+1} = \sum_{i=1}^{n+3}\gra_{i}\E'_i$ and
  $\E_{1} = \sum_{i=1}^{n+3}\grb_{i}\E'_i$. Now \eqref{H01-1} and
  \eqref{H01-3} yield
  \begin{equation*}
    0 \ne \F_1 = \E_{n+1}\E_1 = 
    (\gra_{1}\grb_{2} - \gra_{2}\grb_{1})\F'_3
  \end{equation*}
  and, therefore, $\gra_{1} \ne 0$ or $\gra_{2} \ne 0$. As one has
  \begin{equation*}
    \E_{n+1}\F'_1  = \gra_{1}\E'_{1}\F'_1 = \gra_{1}\G'_1 \qqand 
    \E_{n+1}\F'_2 = \gra_{2}\E'_{2}\F'_2 = \gra_{2}\G'_1
  \end{equation*}
  it follows that $\E_{n+1}\F_j$ is nonzero for some $j$, which
  contradicts the assumption $p=0$; see \eqref{H01-1}.  Thus $\mfB$ is
  not of class $\clB$.

  (g): Assume that $\mfB$ is of class $\mathbf{G}$. The subspace
  \begin{equation*}
    \B_1^\perp \deq \setof{\E \in \B_1}{\E\F = 0 \text{ for all } \F\in\B_2}
  \end{equation*}
  has rank $\prmm{\mfB} - \prmr{\mfB}$.  By \eqref{12} one has
  $\E_{n+2}\F = 0 = \E_{n+3}\F$ for all $\F\in\B_2$, so
  $\B_1^\perp$ has rank at least $2$.
\end{bfhpg}

\begin{bfhpg}[Proof of Proposition \ref{prp:linkH1}]
  \label{prf:linkH1}
  Adopt \stpref{link} and set $p = \prmp{\mfA}$ and $q = \prmq{\mfA}$.
  One may assume that the nonzero products of elements from
  \eqref{Abasis} are
  \begin{equation*}
    \e_1\e_{2+i} = \f_{q+i} \quad\text{for } 1 \le i \le p \qqand
    \e_1\f_j = \g_j \quad\text{for } 1 \le j \le q\:.
  \end{equation*}
  Proceeding as in \pgref{link}, it follows from \pgref{link1} that
  $\mfA$ is directly linked to an ideal $\mfB$ whose Tor algebra
  $\B_\sbt$ has bases
  \begin{equation*}
    \begin{aligned}
      \E_1,\ldots,\E_n,\E_{n+1},\E_{n+3} \quad &\text{for } \B_1\,,\\
      \F_1,\ldots,\F_q,\F_{q+2},\ldots, \F_{l} \quad
      &\text{for } \B_2\,,\text{ and}\\
      \G_4,\ldots,\G_m \quad &\text{for } \B_3\:.
    \end{aligned}
  \end{equation*}
  In particular, one has $\prmm{\mfB} = n + 2$ and
  $\prmn{\mfB} = m-3$.  Further it follows from \eqref{11} and
  \eqref{12} that the nonzero products in $\B_{\sbt}$ that involve
  $\E_{n+1}$ are precisely
  \begin{equation*}
    \tag{1}\label{eq:H1-1}
    \begin{alignedat}{2}
      \E_{n+1}\E_j &= \F_{j}& &\quad\text{for } 1 \le j \le q\quad\text{and}\\
      \E_{n+1}\F_{q+i} &= \G_{2+i}& &\quad\text{for } 2 \le i \le p\:.
    \end{alignedat}
  \end{equation*}
  Set $p' = \prmp{\mfB}$ and $q' = \prmq{\mfB}$; evidently one has
  \begin{equation*}
    \tag{2}\label{eq:H1-2}
    p' \ge q \qqand q' \ge p-1\:.
  \end{equation*}
  Notice that $\mfB$ is not complete intersection, as one has
  $n \ge 2$ by assumption and hence $\prmm{\mfB} \ge 4$; see
  \pgref{ci} and \pgref{gor}. That is, $\mfB$ is of class $\clB$,
  $\mathbf{G}$, $\mathbf{H}$, or $\clT$.\

  (a): If $p \ge 2$ holds, then \eqref{H1-2} yields $q' \ge 1$, and
  per \eqref{pqrmn} it follows that the ideal $\mfB$ is of class
  $\mathbf{B}$, $\mathbf{G}$, or $\mathbf{H}$.

  If $\mfB$ is of class $\clB$, then one has $1 = p' \ge q$, see
  \eqref{pqrmn} and \eqref{H1-2}, so assume towards a contraction that
  $q=0$ holds.  By an argument parallel to the one given in the proof
  of \prpref{linkH0}(a), the nonzero product $\E_{n+1}\F_{q+2} = \G_4$
  from \eqref{H1-1} forces $\E_{n+1}\E_j \ne 0$ for some $j$, which
  contradicts the assumption $q=0$; see~\eqref{H1-1}. Thus
  $q = 1 = p'$ holds.

  If $\mfB$ is of class $\mathbf{G}$, then \eqref{pqrmn} and
  \eqref{H1-2} yield $0 = p' \ge q$, so $p'=0=q$ holds.

  If $\mfB$ is of class $\mathbf{H}$, then an argument parallel to the
  one given in the proof of \prpref{linkH0}(a) shows that the nonzero
  product $\E_{n+1}\F_{q+2} = \G_4$ from \eqref{H1-1} forces $p'$
  linearly independent products of the form $\E_{n+1}\E_j$. This
  implies $p' \le q$, see \eqref{H1-1}, whence $p'=q$ holds by
  \eqref{H1-2}.

  (b): If $p \ge 3$ holds, then \eqref{H1-2} yields $q' \ge 2$, and
  per \eqref{pqrmn} it follows that the ideal $\mfB$ is of class
  $\mathbf{H}$. Moreover, $p'=q$ holds by (a).
  
  (c): Assume that $q \ge 2$ holds. As one has $p'\ge 2$, see
  \eqref{H1-2}, it follows per \eqref{pqrmn} that $\mfB$ is of class
  $\mathbf{H}$ or $\mathbf{T}$. By \eqref{H1-2} it is sufficient to
  prove that $q' \le p-1$ holds.

  If $\mfB$ is of class $\clT$, then per \eqref{pqrmn} one has $0=q'$,
  so $q'\le p-1$ trivially holds.

  If $\mfB$ is of class $\mathbf{H}$, then the argument given in the
  proof of \prpref{linkH0}(c) applies to show that the nontrivial
  products $\E_{n+1}\E_1 = \F_1$ and $\E_{n+1}\E_2 = \F_2$ from
  \eqref{H1-1} force $q'$ linearly independent products of the form
  $\E_{n+1}\F_{j}$. This implies $q'\le p-1$; see \eqref{H1-1}.

  (d): The proof of \prpref{linkH0}(d) applies.

  (e): Assume that $q = 1$ holds. By \eqref{H1-2} one has $p' \ge 1$,
  so the ideal $\mfB$ is per \eqref{pqrmn} of class $\clB$,
  $\mathbf{H}$, or $\clT$.  If $\mfB$ is of class $\clT$, then the
  argument given in the proof of \prpref{linkH0}(e) applies to show
  that the nonzero product $\E_{n+1}\E_1 = \F_1$ from \eqref{H1-1}
  forces two linearly independent products of the form $\E_{n+1}\E_j$,
  which contradicts the assumption $q=1$; see \eqref{H1-1}.

  (f): It follows from (e) that the ideal $\mfB$ is of class $\clB$ or
  $\mathbf{H}$.  If $\mfB$ is of class $\clB$, then the argument in
  the proof of \prpref{linkH0}(f) applies to show that the nonzero
  product $\E_{n+1}\E_1 = \F_1$ from \eqref{H1-1} forces
  $\E_{n+1}\F_j \ne 0$ for some $j$, and that contradicts the
  assumption $p=1$; see \eqref{H1-1}.

  (g): Assume that $n=2$ and $p'=3$ hold. The basis for $\B_1$ is
  $\E_1, \E_2, \E_3, \E_5$, and by \eqref{11} one has $\E_5\E = 0$ for
  all $\E\in\B_1$, so the three products $\E_1\E_2$, $\E_2\E_3$, and
  $\E_3\E_1$ are non-zero. As $\mfB$ is not of class $\clC{3}$, it
  follows that $\mfB$ is of class $\mathbf{T}$; see
  \eqref{Aefg}. Finally, \eqref{H1-1} yields $q=2$ and $p=1$.
\end{bfhpg}

\begin{bfhpg}[Proof of Proposition \ref{prp:linkH2}]
  \label{prf:linkH2}
  Adopt \stpref{link} and set $p = \prmp{\mfA}$ and $q = \prmq{\mfA}$.
  One may assume that the nonzero products of elements from
  \eqref{Abasis} are
  \begin{equation*}
    \e_1\e_{1+i} = \f_{q+i} \quad\text{for } 1 \le i \le p \qqand
    \e_1\f_j = \g_j \quad\text{for } 1 \le j \le q\:.
  \end{equation*}
  Proceeding as in \pgref{link}, it follows from \pgref{link1} that
  $\mfA$ is directly linked to an ideal $\mfB$ whose Tor algebra
  $\B_\sbt$ has bases
  \begin{equation*}
    \begin{aligned}
      \E_1,\ldots,\E_n,\E_{n+1} \quad &\text{for } \B_1\,,\\
      \F_1,\ldots,\F_q,\F_{q+3},\ldots, \F_{l} \quad
      &\text{for } \B_2\,,\text{ and}\\
      \G_4,\ldots,\G_m \quad &\text{for } \B_3\:.
    \end{aligned}
  \end{equation*}
  In particular, one has $\prmm{\mfB} = n + 1$ and
  $\prmn{\mfB} = m-3$.  Further it follows from \eqref{11} and
  \eqref{12} that the nonzero products in $\B_{\sbt}$ that involve
  $\E_{n+1}$ are precisely
  \begin{equation*}
    \tag{1}\label{eq:H2-1}
    \begin{alignedat}{2}
      \E_{n+1}\E_j &= \F_{j}& &\quad\text{for } 1 \le j \le q\quad\text{and}\\
      \E_{n+1}\F_{q+i} &= \G_{1+i}& &\quad\text{for } 3 \le i \le p\:.
    \end{alignedat}
  \end{equation*}
  Set $p' = \prmp{\mfB}$ and $q' = \prmq{\mfB}$; evidently one has
  \begin{equation*}
    \tag{2}\label{eq:H2-2}
    p' \ge q \qqand q' \ge p-2\:.
  \end{equation*}
  If $\prmn{\mfA} = 2$ then $\prmm{\mfB} = 3$ and it follows from
  \pgref{ci} that $\mfB$ is complete intersection. Thus there are
  bases $\set{\E'_1,\E'_2,\E'_3}$ for $\B_1$ and
  $\set{\F'_1,\F'_2,\F'_3}$ for $\B_2$ with non-zero products
  \begin{equation*}
    \E'_{1}\E'_2 = \F'_3\,,\quad \E'_{2}\E'_3 = \F'_1\,,
    \qand \E'_{3}\E'_1 = \F'_2\:.
  \end{equation*}
  For linearly independent vectors
  \begin{equation*}
    \E' \deq \gra_1\E'_1 + \gra_2\E'_2 + \gra_3\E'_3 \qqand
    \E'' \deq \grb_1\E'_1 + \grb_2\E'_2 + \grb_3\E'_3
  \end{equation*}
  the matrix
  \begin{equation*}
    \begin{pmatrix}
      \gra_1 & \gra_2 & \gra_3\\ \grb_1 & \grb_2 & \grb_3
    \end{pmatrix}
  \end{equation*}
  has rank $2$; it follows that the product
  \begin{align*}
    \E'\E'' = (\gra_2\grb_3 - \gra_3\grb_2)\F'_1
    + (\gra_3\grb_1 - \gra_1\grb_3)\F'_2
    + (\gra_1\grb_2 - \gra_2\grb_1)\F'_3
  \end{align*}
  is nonzero. Thus, the products $\E_3\E_1$ and $\E_3\E_2$ are
  nonzero, so \eqref{H2-1} yields $q=2$.

  Assuming now that $m \ge 5$ or $n \ge 3$ and hence
  $\prmm{\mfB} \ge 4$ or $\prmn{\mfB} \ge 2$, it follows \pgref{ci}
  that $\mfB$ is not complete intersection. That is, $\mfB$ is of
  class $\clB$, $\mathbf{G}$, $\mathbf{H}$, or $\clT$.\

  (a): If $p \ge 3$ holds, then \eqref{H2-2} yields $q' \ge 1$ and per
  \eqref{pqrmn} it follows that the ideal $\mfB$ is of class
  $\mathbf{B}$, $\mathbf{G}$, or $\mathbf{H}$.

  If $\mfB$ is of class $\clB$, then one has $1 = p' \ge q$, see
  \eqref{pqrmn} and \eqref{H2-2}, so assume towards a contraction that
  $q=0$ holds.  By an argument parallel to the one given in the proof
  of \prpref{linkH0}(a), the nonzero product $\E_{n+1}\F_{q+3} = \G_4$
  from \eqref{H2-1} forces $\E_{n+1}\E_j \ne 0$ for some $j$, which
  contradicts the assumption $q=0$; see~\eqref{H2-1}. Thus
  $q = 1 = p'$ holds.

  If $\mfB$ is of class $\mathbf{G}$, then \eqref{pqrmn} and
  \eqref{H2-2} yield $0 = p' \ge q$, so $p'=0=q$ holds.

  If $\mfB$ is of class $\mathbf{H}$, then an argument parallel to the
  one given in the proof of \prpref{linkH0}(a) shows that the nonzero
  product $\E_{n+1}\F_{q+3} = \G_4$ from \eqref{H2-1} forces $p'$
  linearly independent products of the form $\E_{n+1}\E_j$. This
  implies $p' \le q$, see \eqref{H2-1}, whence $p'=q$ holds by
  \eqref{H2-2}.

  (b): If $p \ge 4$ holds, then \eqref{H2-2} yields $q' \ge 2$, and
  per \eqref{pqrmn} it follows that the ideal $\mfB$ is of class
  $\mathbf{H}$. Moreover, $p'=q$ holds by (a).
  
  (c): Assume that $q \ge 2$ holds. As one has $p'\ge 2$, see
  \eqref{H2-2}, it follows per \eqref{pqrmn} that $\mfB$ is of class
  $\mathbf{H}$ or $\mathbf{T}$. By \eqref{H2-2} it is sufficient to
  prove that $q' \le p-2$ holds.

  If $\mfB$ is of class $\clT$, then per \eqref{pqrmn} one has $0=q'$,
  so $q'\le p-2$ trivially holds.

  If $\mfB$ is of class $\mathbf{H}$, then the argument given in the
  proof of \prpref{linkH0}(c) applies to show that the nontrivial
  products $\E_{n+1}\E_1 = \F_1$ and $\E_{n+1}\E_2 = \F_2$ from
  \eqref{H2-1} force $q'$ linearly independent products of the form
  $\E_{n+1}\F_{j}$. This implies $q'\le p-2$; see \eqref{H2-1}.

  (d): The proof of \prpref{linkH0}(d) applies.

  (e): Assume that $q = 1$ holds. By \eqref{H2-2} one has $p' \ge 1$,
  so the ideal $\mfB$ is per \eqref{pqrmn} of class $\clB$,
  $\mathbf{H}$, or $\clT$.  If $\mfB$ is of class $\clT$, then the
  argument given in the proof of \prpref{linkH0}(e) applies to show
  that the nonzero product $\E_{n+1}\E_1 = \F_1$ from \eqref{H2-1}
  forces two linearly independent products of the form $\E_{n+1}\E_j$,
  which contradicts the assumption $q=1$; see \eqref{H2-1}.

  (f): It follows from (e) that the ideal $\mfB$ is of class $\clB$ or
  $\mathbf{H}$.  If $\mfB$ is of class $\clB$, then the argument in
  the proof of \prpref{linkH0}(f) applies to show that the nonzero
  product $\E_{n+1}\E_1 = \F_1$ from \eqref{H2-1} forces
  $\E_{n+1}\F_j \ne 0$ for some $j$, and that contradicts the
  assumption $p=2$; see \eqref{H2-1}.
\end{bfhpg}

\bibliographystyle{amsplain-nodash}

\def\soft#1{\leavevmode\setbox0=\hbox{h}\dimen7=\ht0\advance \dimen7
  by-1ex\relax\if t#1\relax\rlap{\raise.6\dimen7
    \hbox{\kern.3ex\char'47}}#1\relax\else\if T#1\relax
  \rlap{\raise.5\dimen7\hbox{\kern1.3ex\char'47}}#1\relax \else\if
  d#1\relax\rlap{\raise.5\dimen7\hbox{\kern.9ex
      \char'47}}#1\relax\else\if D#1\relax\rlap{\raise.5\dimen7
    \hbox{\kern1.4ex\char'47}}#1\relax\else\if l#1\relax
  \rlap{\raise.5\dimen7\hbox{\kern.4ex\char'47}}#1\relax \else\if
  L#1\relax\rlap{\raise.5\dimen7\hbox{\kern.7ex
      \char'47}}#1\relax\else\message{accent \string\soft \space #1
    not defined!}#1\relax\fi\fi\fi\fi\fi\fi}
\providecommand{\MR}[1]{\mbox{\href{http://www.ams.org/mathscinet-getitem?mr=#1}{#1}}}
\renewcommand{\MR}[1]{\mbox{\href{http://www.ams.org/mathscinet-getitem?mr=#1}{#1}}}
\providecommand{\arxiv}[2][AC]{\mbox{\href{http://arxiv.org/abs/#2}{\sf
      arXiv:#2 [math.#1]}}} \def\cprime{$'$}
\providecommand{\bysame}{\leavevmode\hbox to3em{\hrulefill}\thinspace}
\providecommand{\MR}{\relax\ifhmode\unskip\space\fi MR }
\providecommand{\MRhref}[2]{%
  \href{http://www.ams.org/mathscinet-getitem?mr=#1}{#2} }
\providecommand{\href}[2]{#2}

\end{document}